\documentclass[a4paper]{article}
\usepackage{fullpage}
\usepackage{caption}
\DeclareCaptionFont{purple}{\color{purple}} 

\usepackage{amsmath,amsfonts,amsthm,amssymb}
\usepackage{color,xcolor}
\usepackage{ifpdf}
\usepackage{psfrag}
\usepackage{graphicx,graphics,subfigure}
\usepackage{url}
\usepackage{hyperref}
\usepackage{todonotes}

\newtheorem{theorem}{Theorem}[section]

\newtheorem{remark}[theorem]{Remark}

 \newcommand{\N}{\mathbb N}

\newcommand{\modif}[1]{{#1}}
\newcommand{\reviseone}[1]{{#1}}
\newcommand{\revisetwo}[1]{{#1}}

\usepackage{float}
\usepackage{soul}
\usepackage{authblk}
\usepackage[titlenumbered,ruled]{algorithm2e}
\usepackage{comment}
\usepackage{tikz}
\begin{document}
\title{A new variable shape parameter strategy for RBF approximation using neural networks}
\author[1]{Fatemeh Nassajian Mojarrad}
\author[2]{Maria Han Veiga}
\author[3]{Jan S. Hesthaven}
\author[4]{Philipp \"{O}ffner}
\affil[1]{Institute of Mathematics, University of Zürich, Switzerland}
\affil[2]{Department of Mathematics, University of Michigan, USA}
\affil[2]{Michigan Institute for Data Science, University of Michigan, USA}
\affil[3]{Chair of Computational Mathematics and Simulation Science, École Polytechnique Fédérale de Lausanne, Switzerland}
\affil[4]{Institute of Mathematics,  Johannes Gutenberg-University Mainz, Germany}

\date{} 
\maketitle

\begin{abstract}
The choice of the shape parameter highly effects the behaviour of radial basis function (RBF) approximations, as it needs to be selected to balance between \reviseone{the ill-conditioning} of the interpolation matrix and high accuracy. In this paper, we demonstrate how to use neural networks  to determine the shape parameters in RBFs. In particular, we construct a multilayer perceptron \reviseone{(MLP)} trained using an unsupervised learning strategy, and use it to predict shape parameters for inverse multiquadric and Gaussian kernels. We test the neural network approach in RBF interpolation tasks and in a RBF-finite difference method in one and two-space dimensions, demonstrating promising results. 
\end{abstract}

\textbf{keywords.} Meshfree methods; Radial basis function; Artificial neural network; Variable shape parameter; Unsupervised learning

\section{Introduction}

\begin{figure}[ht]
\begin{center}
\includegraphics[width=1.0\textwidth]{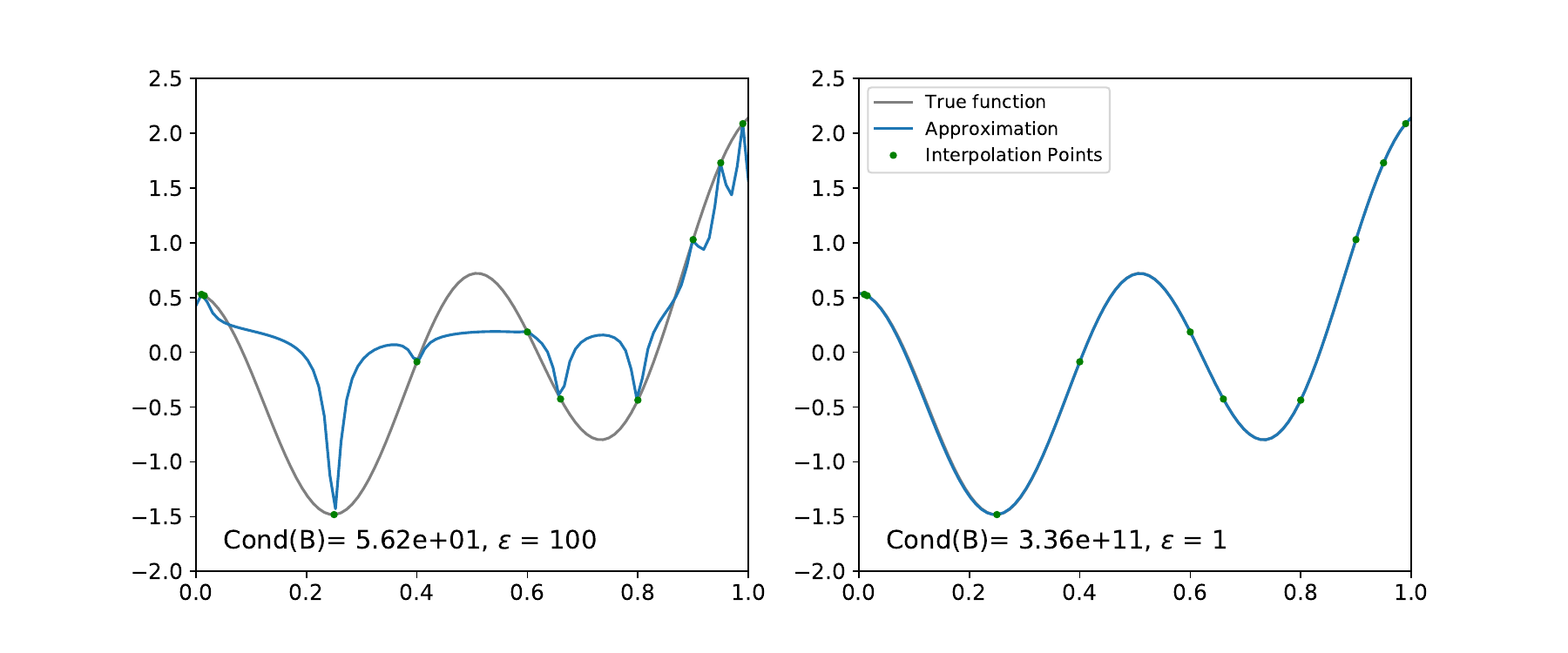}
	\caption{An example of the influence of the choice of $\varepsilon$, and consequently, the condition number of interpolation matrix $B$, using inverse multiquadric basis functions (IMQBF).
	\label{fig:example_rbf}}
\end{center}
\end{figure}

Meshfree methods have great advantages over mesh-based methods, since they do not require a grid based domain or surface discretization. The RBF-interpolation is known as a truly meshfree computational method and has the advantage that it is simple to implement and can be applied in complex geometries in high dimensions. Since the \reviseone{RBF values} depend only on the distance to its center point, RBFs can be easily implemented to reconstruct a hypersurface using scattered data in multi-dimensional spaces \cite{Fasshauer, Flyer, wendland2004scattered, lazzaro2002radial}. Additionally, for smooth functions, spectral convergence \reviseone{can be achieved} \cite{fornberg2005accuracy, flyer2016role}.\\
Globally supported RBFs often \reviseone{include} a shape parameter that plays an important role in the accuracy of the approximation. However, it remains a challenge to find the optimal shape parameter \reviseone{value} \cite{O.Davydov,B.Fornberg,M.Mongillo,Ranjbar}. \reviseone{As first noted in \cite{schaback1995}, there exists a \textit{trade-off} between the error and the condition number of the interpolation matrix (sensitivity), in the sense that there is no known case where both the interpolation error and the interpolation matrix condition number are small. Either one optimizes for a small error and gets a near ill-conditioned matrix, or one chooses an interpolation matrix with a small matrix condition number, at the cost of having a larger interpolation error.} In contrast, a large shape parameter results in well-conditioned matrices, but adversely impacts the accuracy. A depiction of the influence of the choice on the shape parameter $\varepsilon$ on the approximation property and condition number of the interpolation matrix is shown in Fig. \ref{fig:example_rbf}. Thus, the choice of the shape parameter plays an important role and affects both the interpolation error and the stability condition of the RBF. \revisetwo{In the context of PDEs, for example, in \cite{C.Piret} a study on how the accuracy varies with RBF type, shape parameter and length of time integration is provided when solving convective PDEs.}

There \reviseone{have been many} attempts to address this problem. For example, one of the early adaptive approaches was proposed by Hardy \cite{hardyrbf}, suggesting $\varepsilon = 1/(0.815d)$, with $d=1/N \sum_{i=1}^N d_i $, $N$ being the total number of interpolation nodes and $d_i$ being the distance of datapoint $i$ to its nearest neighbour. Another early proposal for adaptive $\varepsilon$ was by Franke \cite{frankerbf}, giving $\varepsilon = \frac{0.8\sqrt{N}}{D}$, where $D$ is the diameter of the smallest circle containing all data points.

\revisetwo{Other approaches rely on locally adjusting the shape parameters. For example, in \cite{M.Mongillo}, predictor functions are introduced to mimic the interpolation error, and to search for an optimal shape parameter $\varepsilon$, whereas in \cite{O.Davydov}, a multilevel algorithm that finds a near-optimal shape parameter is presented. Another way to \revisetwo{address this} problem is by evaluating the shape parameter expansion with more advanced algorithms, e.g. allowing complex values or using an efficient rational approximation technique  \cite{B.Fornberg,G.B.Wright}. }



\revisetwo{However, even if these ansatzes show promising results, disadvantages remain, e.g. they are not computationally fast nor always very accurate, and for any set of data, the strategy has to be re-applied.} In this paper we present a new \revisetwo{simpler and general} approach  to select the shape parameter via neural networks (NN), using an unsupervised learning strategy. Our strategy is to learn a map that predicts a shape parameter $\varepsilon$. \revisetwo{In particular, contrasting with the previously mentioned works, once the NN has been trained successfully, the evaluation of the NN is extremely quick, leading to negligible overhead in comparison to a fixed shape parameter, as shown in our numerical experiments.}

The outline of the manuscript is as follows: in Section \ref{se_RBF}, we briefly review the RBF method for interpolation and the \reviseone{RBF} based finite difference (RBF-FD) method for both steady-state and time-dependent problems. In Section \ref{se_MLP}, we introduce some basic concepts of \reviseone{NNs} and demonstrate the use of \reviseone{NNs} to predict the optimal local values of the shape parameter to control the accuracy. In \reviseone{Section} \ref{se_results} we test the validity and effectiveness of the proposed strategy, namely, we will test this idea on interpolation problems and solving partial differential equations \reviseone{(PDEs)}. 
Finally, \reviseone{Section} \ref{se_discussion} contains a discussion about the obtained results, as well as future work, and  provides some conclusions.

\section{Radial Basis Functions (RBF)}\label{se_RBF}

RBFs can be used for a variety of tasks. They are good alternatives to interpolate functions on arbitrary sets of multi-dimensional points which, due to the Mairhuber-Curtis theorem, \reviseone{is not} possible with (multivariate) polynomials of degree $N$ in $\mathbb{R}^d$, $d>1$ \cite{zbMATH03117700,curties_1}.

\subsection{RBF interpolation}
In this section, we introduce the RBF method for the interpolation of scattered data. 
For a set of N distinct centers $\mathbf{x}_1,\mathbf{x}_2,\dots,\mathbf{x}_N$ in $\Omega\subset \mathbb{R}^d$, the RBF interpolation function  can be written in the form
\begin{equation}
\label{eq:kansa}
  S(\textbf{x})=\sum_{i=1}^N \lambda_{i}\phi(||\mathbf{x}-\mathbf{x}_i||) +\sum_{k=1}^m \gamma_{k}\bar{p}_k(\mathbf{x}), 
\end{equation}
where $||.||$ denotes the standard Euclidean norm, $\phi$ is the \reviseone{RBF} (also called kernel), $ \lambda_{i}$ represents the unknown  interpolation coefficients, $\gamma_k$ the unknown Lagrange
multipliers and $\{ \bar{p}_k\}_{k=1}^m$ is a basis for the space of polynomials up to degree $m-1$, denoted by $\mathcal{P}_{m-1}$.
The number of monomial terms is \modif{$m = \binom  {D_m+(d-1)}{d-1}$, where $D_m$ is the degree of the monomial basis and $d$ the dimension.}
The RBF interpolant \eqref{eq:kansa} of a function $f$ is uniquely determined by the conditions 
\begin{equation}\label{eq_condition}
    \begin{split}
     &S(\textbf{x}_i)=f(\textbf{x}_i), \qquad i=1, \cdots, N, \\
     &\sum_{i=1}^N \lambda_{i} \bar{p}_k(\mathbf{x}_i)=0, \qquad k=1,\cdots,m. 
    \end{split}
\end{equation}
which are used to calculate the unknown coefficients $\mathbf{\lambda}=(\lambda_{1}, \cdots, \lambda_N)^T$ and $\mathbf{\gamma}= (\gamma_{1}, \cdots, \gamma_m)^T$. 
We can reformulate \eqref{eq_condition} as a linear system to calculate the unknowns, \reviseone{namely,}
\begin{equation}
\label{eq:approximation_prob}
 \begin{bmatrix}
 A &P\\
 P^T&0
 \end{bmatrix}
 \begin{bmatrix}
 \lambda\\
 \gamma
 \end{bmatrix}=
 \begin{bmatrix}
 \mathbf{f}\\
 0
 \end{bmatrix},
\end{equation}
where $\mathbf{f}=(f(\textbf{x}_1), \cdots,f(\textbf{x}_N ))^T$ and 
\begin{equation}
\label{eq:interpolation matrix}
A=
  \begin{bmatrix}
  \phi(||\mathbf{x}_1-\mathbf{x}_1||) &   \phi(||\mathbf{x}_1-\mathbf{x}_2||)&  \dots &  \phi(||\mathbf{x}_1-\mathbf{x}_N||)\\
      \phi(||\mathbf{x}_2-\mathbf{x}_1||) &   \phi(||\mathbf{x}_2-\mathbf{x}_2||)&  \dots &  \phi(||\mathbf{x}_2-\mathbf{x}_N||)\\
    \dots  &  \dots& \dots &  \dots \\
 \phi(||\mathbf{x}_N-\mathbf{x}_1||) &   \phi(||\mathbf{x}_N-\mathbf{x}_2||)&  \dots &  \phi(||\mathbf{x}_N-\mathbf{x}_N||)    
  \end{bmatrix}
,\quad 
P=
  \begin{bmatrix}
  \bar{p}_1(\textbf{x}_1)&  \dots & \bar{p}_m(\textbf{x}_1)\\
 \bar{p}_1(\textbf{x}_2)&  \dots & \bar{p}_m(\textbf{x}_2) \\
 \dots &\dots&\dots\\
\bar{p}_1(\textbf{x}_N) &\dots & \bar{p}_m(\textbf{x}_N)  
  \end{bmatrix}.
\end{equation}
We denote by $B=\begin{bmatrix}
 A &P\\
 P^T&0
 \end{bmatrix}$ and assume that the solution \reviseone{at} any point $\mathbf{y}\in \Omega$ can be written in the form
 \begin{equation}
 \begin{split}
  S(\textbf{y} )&= \underbrace{(\phi(||\mathbf{y}-\mathbf{x}_1||),\cdots,\phi(||\mathbf{y}-\mathbf{x}_N||),\bar{p}_1(\mathbf{y}),\cdots,\bar{p}_m(\mathbf{y}))}_{q(\mathbf{y})} \begin{bmatrix}
 \lambda\\
 \gamma
 \end{bmatrix}\\
 &=(q(\mathbf{y})B^{-1})_{1:N}\mathbf{f}=\sum_{i=1}^N f(\textbf{x}_i)[q(\mathbf{y})B^{-1}]_i =\sum_{i=1}^N f(\textbf{x}_i) \psi_i(\mathbf{y}),
 \end{split}
 \label{eq:eq1}
 \end{equation}
 where $ \psi_i$ are called  cardinal functions. They \reviseone{form} a basis of the underlying approximation space.\\
Commonly used types of \reviseone{RBFs} are listed in Table  \ref{tab:1}. 
In this work, we use the inverse multiquadric RBF
$$\phi(r)=\frac{1}{\sqrt{1+(\varepsilon r)^2}},$$
 and the Gaussian RBF \reviseone{(GARBF)} with shape parameter $\varepsilon>0$. The condition number of matrix $B$ depends on the
shape parameter $\varepsilon$ as illustrated in Fig. \ref{fig:example_rbf}.
\begin{table}[H]
\centering
\caption{Common RBF types.}
\begin{tabular}{|c| c |} 
\hline
Type of basis function&$\phi(r)$\\
\hline
Gaussian &$e^{-(\varepsilon r)^2}$ \\

Multiquadric &$\sqrt{1+(\varepsilon r)^2}$ \\

Inverse quadratic &$\frac{1}{1+(\varepsilon r)^2}$ \\

Inverse multiquadric &$\frac{1}{\sqrt{1+(\varepsilon r)^2}}$ \\

Polyharmonic spline &$ r^k,~k=1,3,5,\dots$ \\
 &$ r^k \ln(r),~k=2,4,6,\dots$ \\
\hline
\end{tabular}
\label{tab:1}
\end{table}
\begin{remark}
Inclusion of polynomials up to a certain degree inside the RBF interpolation \eqref{eq:kansa} has several advantages: 
\begin{itemize}
    \item The RBF interpolants become exact for polynomials up to a certain degree. Including constants is necessary for the RBF method to be conservative \reviseone{in context of hyperbolic conservation laws} \cite{glaubitz2022energy, glaubitz2022summation}. 
    \item For some (unconditionally) positive kernels $\phi$, the RBF interpolant exists uniquely when polynomials up to a certain degree are included \cite{Fasshauer}.
    \item It has a positive influence on the accuracy of the RBF-FD method \cite{flyer2016role}. 
\end{itemize}
In our work, we will include constant polynomials in our RBF  interpolant \eqref{eq:kansa}. 
\end{remark}

\subsection{RBF-FD method for solving steady-state and time dependent differential equations}

The RBF-FD method is a meshfree method to solve \reviseone{PDEs} numerically, which arose from the global RBF methods \cite{Flyer}. Solving PDEs with global RBF methods results in remarkable numerical aspects like time stability and accuracy. However, the high computational cost and memory requirements of \reviseone{RBFs} for very large problem sizes is the main reason for the development of \reviseone{the} RBF-FD method \cite{Flyer}. RBF-FD formulas are generated from RBF interpolation over local sets of nodes on the surface \cite{Tolstykh}, i.e.: only a nearest neighbour region is connected to each basis function. This means that we  have  a clustering of point clouds where the cardinal basis are considered locally inside the cloud.  This leads to a sparse matrix instead of the full linear system of equations. Next, we derive the RBF-FD formulation \cite{Mishra} for Poisson and time dependent differential equations. 

We are interested in the solution of PDEs of the form:
\begin{equation}
\label{eq:poisson}
  \Delta u(\mathbf{x})=f(\mathbf{x}) \quad  \text{ or} \quad \mathcal{L} u(\mathbf{x})=f(\mathbf{x}),~ \mathbf{x} \in \Omega,
\end{equation}
where  $\Delta $ is the Laplace operator and $\mathcal{L}$ denotes the space differential operator which can be also used in a more general setting.
Equation \eqref{eq:poisson} is complemented by boundary conditions. 
The domain $\Omega$ is discretized using two \reviseone{sets of points}:
\begin{itemize}
\item The interpolation point set $X=\{\textbf{x}_i\}_{i=1}^N$ for generating the cardinal functions.
\item The evaluation point set $Y=\{\textbf{y}_j\}_{j=1}^M$ for sampling the PDE \eqref{eq:poisson}.
\end{itemize}
\reviseone{The difference between the two sets is that $X$ is used to build the underlying approximation space whereas at the points $Y$, the PDE is actually solved. If $X=Y$ one speaks about a collocation method and when $M>N$, it is referred as oversampling. }
The cardinality of \reviseone{the} $X$ and $Y$ point sets  matches the relation $M=rN$, where $r$ is the oversampling parameter. We use collocation and oversampled discretizations. 
\reviseone{To describe the algorithm to solve the PDE, we start using the nodal points $X$. We focus first on the time independent  problem for simplicity.
We form a stencil for each $x_i \in X$, based on the $n$ nearest neighbors. Each stencil (also referred as node's support domain)  is represented with a point set $\mathcal{N}(\mathbf{x}_i) \subset X, i=1, \cdots N.$ }
\reviseone{On each stencil point set $\mathcal{N}(\mathbf{x}_i)$, we form a local interpolation problem leading to  
a set of stencil based approximations given by 
}
\begin{equation}\label{eq_local}
 u^{(i)}(\textbf{y})=\sum_{k=1}^n u(\textbf{x}^{(i)}_k)[q(\mathbf{y},\mathcal{N}(\mathbf{x}_i))B(\mathcal{N}(\mathbf{x}_i),\mathcal{N}(\mathbf{x}_i))^{-1}]_k=\sum_{k=1}^n u(\textbf{x}^{(i)}_k) \psi^{(i)}_k(\mathbf{y}), \quad \forall i=1, \cdots, N
\end{equation}
\reviseone{where $\psi^{(i)}_k(\mathbf{y})$
are the local  RBF-FD cardinal functions and $u^{(i)}(\mathbf{y})$ is the local approximation. Note that \eqref{eq_local} is analogous to  \eqref{eq:eq1} where $N$ is replaced by $n$, $f$ by $u$ and $X$ by $\mathcal{N}(\mathbf{x}_i)$.  }

To obtain a global solution at $\mathbf{y}$, we associate every evaluation point  with an index of the closest stencil center point defined by 
\begin{equation*}
  \nu(\revisetwo{\textbf{y}})=\mathop{\arg \min}\limits_{i\in\{1,\cdots,N\}} ||\textbf{y} -\textbf{x}_i|| 
\end{equation*}
This means that $\textbf{y} \in \mathcal{N}(\mathbf{x}_{\nu(\textbf{y})})$, and we have
\begin{equation}\label{eq_local_2}
 u(\textbf{y})=u^{(\nu(\textbf{y}))}(\textbf{y})=\sum_{k=1}^n u(\textbf{x}^{(\nu(\textbf{y}))}_k) \psi^{(\nu(\textbf{y}))}_k(\mathbf{y}).   
\end{equation}
\reviseone{This means nothing else that the global solution is represented by the local solutions inside each region\footnote{We have around each center point $\mathbf{x}_i$ a dual Voronoi diagram hidden, which assigns each $\mathbf{y}$ to a specific local approximation.}. Finally, to get a global representation, we may introduce an auxiliary   index transformation $\sigma$ to rewrite  \eqref{eq_local_2} in the following way }
\begin{equation}\label{eq_global_RBF}
 u(\textbf{y})=\sum_{k=1}^n u(\textbf{x}^{(\nu(\textbf{y}))}_k) \psi^{(\nu(\textbf{y}))}_k(\mathbf{y})=
 \sum_{i=1}^N u(\textbf{x}_i) \psi^{(\nu(\textbf{y}))}_{\sigma(\nu(\mathbf{y}),i)}(\mathbf{y})=\sum_{i=1}^N u(\textbf{x}_i) \psi_i(\mathbf{y})
\end{equation}
\reviseone{where $\psi_i=\psi^{(\nu(\textbf{y}))}_{\sigma(\nu(\mathbf{y}),i)}$ denote the global RBF-FD cardinal functions. Here, $\sigma$ gives the local index in the corresponding point set.}
It follows in a similar fashion that 
\begin{equation}\label{eq:laplace_operator}
 \mathcal{L}u(\textbf{y})= \mathcal{L} u^{(\nu(\textbf{y}))}(\textbf{y})=\sum_{k=1}^n u(\textbf{x}^{(\nu(\textbf{y}))}_k)  \mathcal{L}\psi^{(\nu(\textbf{y}))}_k(\mathbf{y}) = \sum_{i=1}^N u(\textbf{x}_i)  \mathcal{L}\psi_i(\mathbf{y}),   
\end{equation}
which is used in the discretization. 
\reviseone{This means in the steady-state case, we solve the system of equations on the evaluation points using the following discretization:
\begin{equation*}
\sum_{i=1}^N u(\textbf{x}_i)  \mathcal{L}\psi_i(\mathbf{y}_j)=f(\mathbf{y}_j), \quad  \forall j=1, \cdots M. 
\end{equation*}
}

Let us now consider the linear \reviseone{initial-boundary} value problem 
\begin{align}
\begin{split}\label{n2}
\frac{\partial u(\mathbf{x},t)}{\partial t}= {}&\mathcal{L}u(\mathbf{x},t),~ \mathbf{x} \in \Omega,  ~0\leq t\leq T,
\end{split}\\
\begin{split}\label{n5}
u(\mathbf{x},t)={}& g(\mathbf{x},t),~ \mathbf{x} \in \partial \Omega ,
\end{split}\\
u(\mathbf{x},0)={}& f(\mathbf{x}),~\mathbf{x}  \in \Omega  \label{n6}.
\end{align}
where $\Omega \subset  \mathbb{R}^d$ is the solution domain with the boundary $\partial \Omega$ where Dirichlet boundary conditions are considered and $\mathcal{L}$ is a differential operator on $\mathbb{R}$.\\
\reviseone{The difference from before is the dependence of $u$ on space and time now.}
By using a simple method of lines (MOL) approach, we split space and time discretization \reviseone{in the usual way.} For the space discretization our local RRB-FD discretization is applied \reviseone{where the coefficients depend  on time as well and are updated via } classical time-integration schemes like Runge-Kutta (RK) methods, backward differentiation formulas (BDFs) or deferred correction schemes (DeC). 
Analogous to \eqref{eq:laplace_operator}, we have now 
\begin{equation}\label{eq:diff_operator}
 \mathcal{L}u(\textbf{y},t)= \mathcal{L}u^{(\nu(\textbf{y}))}(\textbf{y},t)=\sum_{k=1}^n u(\textbf{x}^{(\nu(\textbf{y}))}_k,t)  \mathcal{L}\psi^{(\nu(\textbf{y}))}_k(\mathbf{y}) = 
 \reviseone{\sum_{i=1}^N u(\textbf{x}_i,t)  \mathcal{L}\psi_i(\mathbf{y})}
\end{equation}
and we obtain a linear system of ordinary \reviseone{differential} equations from \eqref{n2} \reviseone{that involve} the time dependent coefficients $u(\textbf{x}^{(\nu(\textbf{y}))}_k,t)$ \reviseone{(or $u(\textbf{x}_i,t))$}.
In this work we consider \eqref{eq:poisson} or \eqref{n2} with  Dirichlet boundary conditions. We impose the boundary conditions in
two ways: (a) when studying the eigenvalue spectra, we remove the values of the vector $u(\textbf{x}, t)$ located on the boundary,
along with the corresponding columns and rows in the matrices, (b) when computing the numerical solution in time, we use the injection
method to impose the condition, i.e., we overwrite the solution $u(\textbf{x}, t_i), i = 0, 1, \cdots$ with the corresponding
 boundary value, after each step of the  time-stepping algorithm. Note that these two
procedures  are equivalent.

\section{Methodology}\label{se_MLP}
\reviseone{As denoted before, there exists a trade-off between the interpolation error and the interpolation matrix condition number, namely, that the interpolation error appears to decrease with the decrease of the sensitivity of the interpolation matrix, given by a growing matrix condition number.} 

\reviseone{We want to find a mapping between a set of nodal  points and an optimal shape parameter $\varepsilon$ that satisfies accuracy and stability of the system. To achieve this, we use a data-driven approach.}

Therefore, we formulate our \reviseone{learning} goal as follows: find a function $f$ for which given a set of $n$ arbitrary points, denoted as a vector $\vec{x}\in\mathbb{R}^{n}$, returns $\varepsilon$ such that the condition number of the interpolation matrix $M(\vec{x},\varepsilon)$, \reviseone{ which is generated} using radial basis functions \reviseone{with parameter $\varepsilon$} and set of points $\vec{x}$, is close to some large \reviseone{(but acceptable)} value $\text{cond}_{\text{max}}$ at the \textit{edge of} being ill-conditioned:

\begin{equation}
    \label{eq:approx function}
    f:\vec{x}\in \mathbb{R}^n \to \varepsilon \quad \text{s.t. } \text{cond}(M(\vec{x},\varepsilon))\approx  \text{cond}_{\max}.
\end{equation}

\reviseone{In the rest of this section, we will describe the dataset generation (Section \ref{sec:dataset_generation}), the data-driven model approximating \eqref{eq:approx function} and the integration of the method with the aforementioned tasks, namely, function interpolation and RBF-FD.}

\subsection{Dataset generation}
\label{sec:dataset_generation}
\reviseone{Focusing on the 1-dimensional case,} we consider $m$ sets of $n$ arbitrary points in the interval $[0,1]$. These define our interpolation nodes, which we will refer to as \textit{coordinate-based features}\reviseone{, yielding a set of feature vectors $\{\tilde{x}^i=(x_1 ^i,...,x_n ^i)\}_{i=1}^m$.}

We can also generate a different set of features, called \textit{distance-based features}, which is generated from the interpolation nodes. We define $n-1$ distances as:
\begin{align*}
        &d_1 ^i= x_2^i - x_1 ^i\\
        &\cdots  \\
        &d_{n-1}^i = x_{n}^i - x_{n-1}^i, \\
\end{align*}
We note that $x_p ^i-x_q ^i = \sum_{l=q}^{p-1} d_l ^i$ \revisetwo{and for $p=q$, we have always $d_0 ^i = 0$}. Thus, \reviseone{we} can \reviseone{express} the interpolation matrix \eqref{eq:interpolation matrix} in terms of the distance feature vector $x^i=(d_1^i,...,d_{n-1}^i)$.

For the 1-dimensional problem, it is advantageous to use distance-based features as the dimension is reduced and we have a shift-invariance property. We also observe experimentally that the distance-based features appear to lead to a better performance across different problems.

\reviseone{Finally, we standardise the data by considering the inverse of the distances, and then centering and normalising all data points:}
\begin{equation}
\label{eq:transform}
x_j'^i = \frac{x_j^i-\bar{x}_j}{{\sigma_j}^2}, \quad i = 1,...,m, \quad j=1,...,k
\end{equation}
for \reviseone{$k$} denoting the number of features, $\bar{x}_j$ the empirical mean and ${\sigma_j}^2$ the empirical variance across all training points.

For \reviseone{the} construction of a training sample, we consider samples of $n$ uniformly distributed points over $[0,1]$. 
The training and validation datasets are constructed using the distances between these points. We consider a total of 4400 training samples and 1100 validation samples.

\subsection{Data driven model}
\label{sec:model_data_driven}

\def\layersep{3cm}
\def\nodeinlayersep{1.5cm}
\begin{figure}
\centering
\begin{tikzpicture}[
   shorten >=1pt,->,
   draw=black!50,
    node distance=\layersep,
    every pin edge/.style={<-,shorten <=1pt},
    neuron/.style={circle,fill=black!25,minimum size=17pt,inner sep=0pt},
    input neuron/.style={neuron, fill=green!50},
    output neuron/.style={neuron, fill=red!50},
    hidden neuron/.style={neuron, fill=blue!50},
    annot/.style={text width=4em, text centered}
]
    \foreach \name / \y in {1,...,3}{
    \ifnum \y=2
        \node at (0,-\y-2.5) {$\vdots$};
    \else
        \node[input neuron] (I-\name) at (0,-\y-2.5) {};  
    \fi
    }
    \newcommand\Nhidden{2}

    \foreach \N in {1,...,\Nhidden} {
       \foreach \y in {1,...,5} { 
         \ifnum \y=3
           \node (dots\N-\y) at (\N*\layersep,-\y*\nodeinlayersep) {$\vdots$};
         \else
           \node[hidden neuron] (H\N-\y) at (\N*\layersep,-\y*\nodeinlayersep ) {};
         \fi
       }
    \node[annot,above of=H\N-1, node distance=1cm] (hl\N) {Hidden layer \N};
    }
    \node[output neuron, right of=dots\Nhidden-3] (O) {$\varepsilon$}; 
    \foreach \source in {1,3}
        \foreach \dest in {1,2,4,5} 
            \path (I-\source) edge (H1-\dest);
    \foreach [remember=\N as \lastN (initially 1)] \N in {2,...,\Nhidden}
       \foreach \source in {1,...,2,4,5} 
           \foreach \dest in {1,...,2,4,5} 
               \path (H\lastN-\source) edge (H\N-\dest);
    \foreach \source in {1,...,2,4,5} 
        \path (H\Nhidden-\source) edge (O);
    \node[annot,left of=hl1] {Input layer};
    \node[annot,right of=hl\Nhidden] {Output layer};
\end{tikzpicture}

\caption{\revisetwo{An example of a dense neural network, with two hidden layers (i.e.: $p-1=2$). The width of the input layer depends on the choice of features, namely in the 1-dimensional setting, $q_0=n$ for coordinate based features and $q_0=n-1$ for distance based features. The output layer returns the shape parameter $\varepsilon$.}\label{fig:nn_diagram}}
\end{figure}
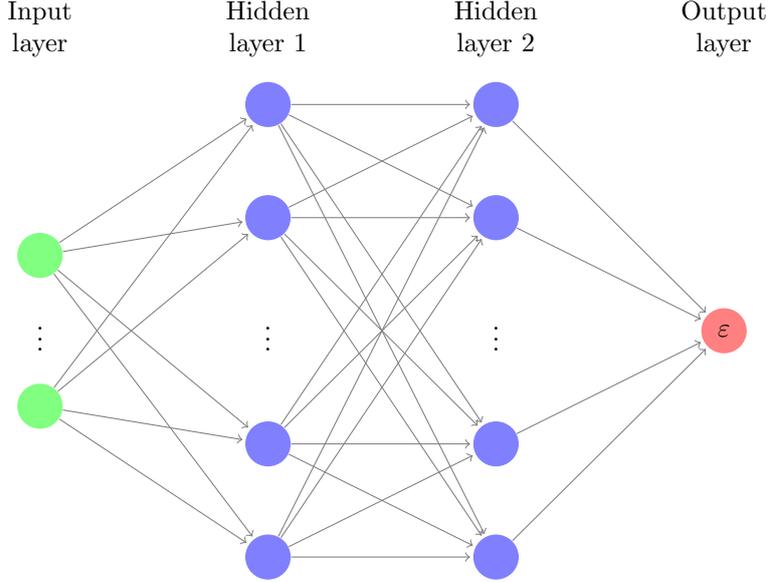

\reviseone{To approximate the map in \eqref{eq:approx function}, we consider a dense, fully connected NN}. We assume that $f$ is approximated by function $\mathcal{H}$ of the form:
\begin{equation}
 \mathcal{H}(\vec{x}) := \sigma_{p}(...\sigma_1(A_1 \vec{x} + b_1)),
\end{equation}
for activation functions $\sigma_i$, and trainable parameters $\mathbf{W}=\{A_i,b_i\}$, for $i=1,\cdots,p$. $A_i$ is a matrix in \revisetwo{$\mathbb{R}^{q_i\times q_{i-1}}$ and $b_i$ is vector in $\mathbb{R}^{q_i}$, $q_0$ is the input dimension (i.e. $n$ for coordinate-based features and $n-1$ for distance-based features), and $d_p$ is the output dimension ($q_p=1$).} \revisetwo{We show in Figure \ref{fig:nn_diagram} for an example of a plausible NN.} \reviseone{The activation function plays an
important role in enhancing the networks learning capabilities. A popular activation function used by most practitioners, is the
rectified linear unit (ReLU):
\begin{equation*}
    \sigma(x)=\max(0,x).
\end{equation*}
When compared to models using the logistic or hyperbolic tangent activation function, training
with ReLU is often faster, for reasons of linearity, non-saturating form of the ReLU function, and the inexpensive evaluation of the function.}

The structure of the \reviseone{NN} for $n=10$ and  distance-based features has the following parameters: $p=4$, $d_0=9$ , $d_1=14$, $d_2=8$, $d_3=3$, $d_4=1$, with ReLU activation function for $p=1,2,3$ and linear activation for $p=4$.

\reviseone{To encode the objective described in the beginning of the section, we design an appropriate loss function.} We train the \reviseone{NN} with respect to the following cost function, for a data point $\vec{x}_j$:
\begin{equation}
    \label{eq: cost function}
    C_j(\mathbf{W})=
\begin{cases}
0.1 \times\log\big(-\text{cond}(B_j)+10^{10}+\kappa\big),&\text{if } \text{cond}(B_j)\leq 10^{10}, \\
0, &\text{if } 10^{10}<\text{cond}(B_j)\leq10^{12}, \\
\log\big(\frac{10}{9}(\text{cond}(B_j)-10^{12})+\kappa\big), &\text{if } 10^{12}<\text{cond}(B_j)\leq10^{13},\\
 \log\big(\text{cond}(B_j)+\kappa\big),&\text{if } \text{cond}(B_j)> 10^{13},
\end{cases}
\end{equation}
where $\kappa=1$ and $\text{cond}(B_j)=||B_j||_F~ ||B^{-1}_j||_F$, $B_j$ is the interpolation matrix built using the interpolation nodes $\vec{x}_j$ and the predicted $\varepsilon$. \reviseone{We ensure that the loss function is continuous and that privileges $\varepsilon$ that lead to a condition number $B_j$ that is between $10^{10}$ and $10^{12}$.}  Lastly, since we work with condition numbers and we use a gradient descent based method to compute the model parameters, we consider a logarithmic cost function to have a smaller magnitude update in the weights at each iteration of the optimisation algorithm.

Furthermore, special care must be taken to ensure that the trained network does not overfit the training dataset. A common method to avoid overfitting is to regularize the cost functional by penalizing the parameters $\mathbf{W}$ of the network \cite{Nowlan}
\begin{equation}
  \mathcal{R}(\mathbf{W}) = \beta ~||\mathbf{W}||_2^2,
  \label{regularization}
\end{equation}
where $\beta \geq0 $ is the regularization parameter and $||\mathbf{W}||_2^2$ represents the squared sum of all the weights in the MLP.

\reviseone{
Thus, the full cost function is given by:}

\begin{equation}
\label{eq:fullcost}
    C(\mathbf{W})=\frac{1}{m}\sum_{j=1}^{m} C_j(\mathbf{W}) + \beta ~||\mathbf{W}||_2 ^2.
\end{equation}
\reviseone{
We minimize the cost function \eqref{eq:fullcost} with respect to the parameters $\mathbf{W}$ of the NN using stochastic gradient descent, considering a batch size of $500$.}

\begin{remark}
We emphasize that this learning task does not fall into the standard supervised learning framework. We do not have access to an optimal $\varepsilon$ for each set of nodes $\vec{x}$. The network is optimised with respect to the defined loss function by finding a suitable $\varepsilon$ that leads to a matrix with a desired matrix condition number. We set this condition number to be close to $10^{12}$.
\end{remark}

\begin{figure}[H]
\begin{center}
\subfigure[train loss]
{\includegraphics[width=0.45\textwidth]{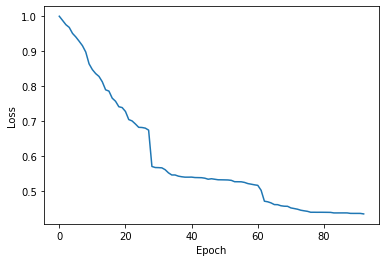}}
\subfigure[validation loss]
{\includegraphics[width=0.45\textwidth]{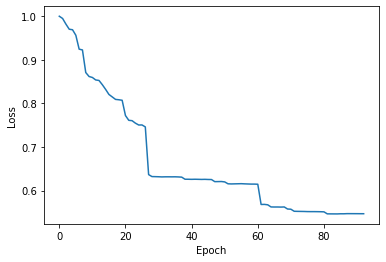}}
\end{center}
	\caption{Loss during training in 1D with IMQRBF.
}
	\label{fig:train_valid_loss}
\end{figure}

In addition, we have an  early stopping \reviseone{criterion} \cite{Nowlan}. We use a validation data set which is independent of the training set. Once the loss computed on the validation data set stagnates after $M$ consecutive epochs, the training is terminated and we recover the best set of parameters. Fig. \ref{fig:train_valid_loss} shows the plots of training loss and validation loss. 
\subsection{Extension to 2-dimensional problems}
We also consider an extension of this method to 2-dimensional problems. For the dataset generation and representation, we have again two options: \textit{coordinate-based features} and \textit{distance-based features}.

In a 2-dimensional problem, the advantage of considering coordinate-based features is clearer, because the number of features per number of interpolation nodes scales as $2\times n$, whereas for distances, we have $\binom{n}{2}$ features. However, when adopting a coordinate-based feature representation we lose the shift-invariant property.

\begin{figure}[ht]
\begin{center}
{\includegraphics[width=0.5\textwidth]{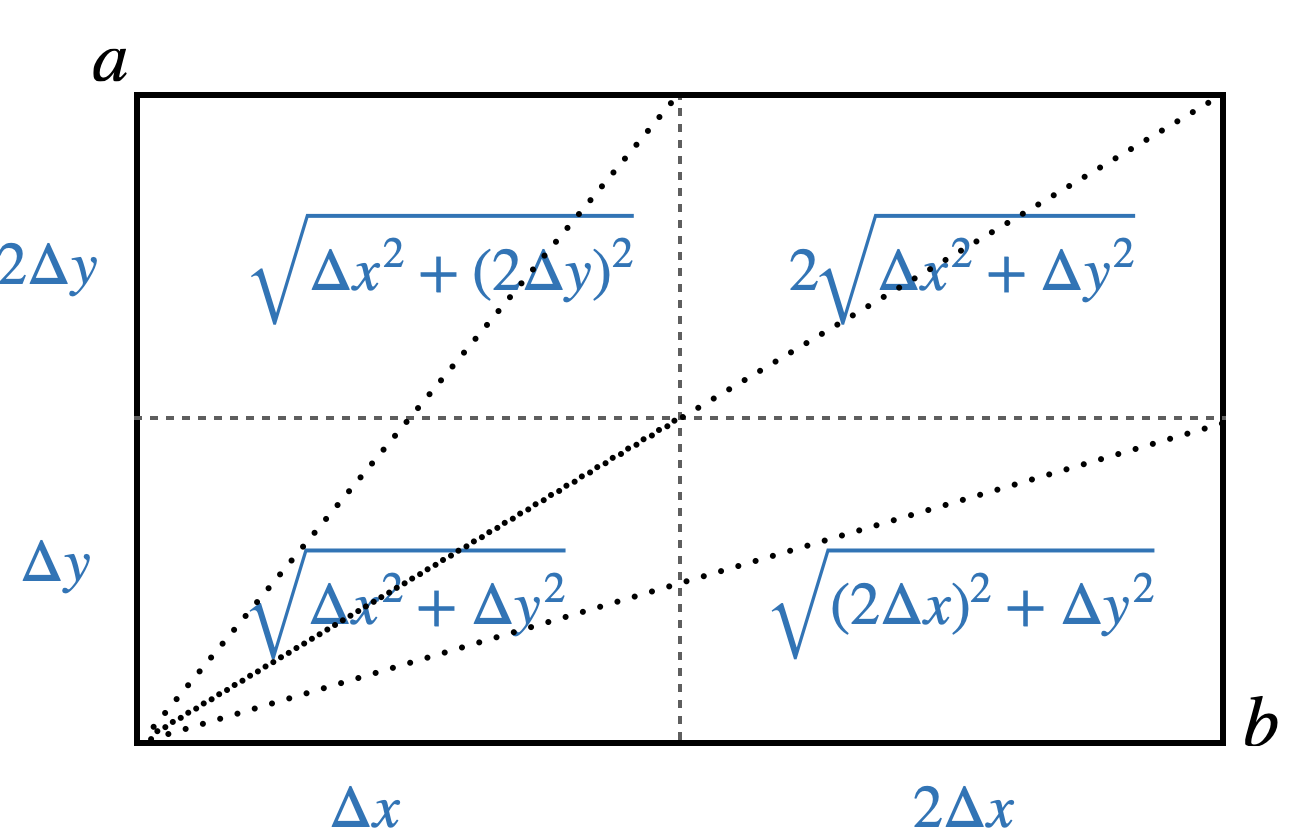}}
	\caption{Structured distribution for 9 nodes. The parameters $a, b \in [0,1]$ define the shape of the rectangle, $\Delta x=b/2$, $\Delta y=a/2$.}
	\label{fig: 2d points}
\end{center}
\end{figure}

In our work, we extend the proposed methodology to a non-generic distribution of 2-dimensional points. In particular, we consider only structured distributions of points (as shown in Fig. \ref{fig: 2d points}, for $n=9$). This constraint on the distribution of points allows the simplification of the description of the cloud of points when considering distance-based features. 

The \revisetwo{feature vector for the 2-dimensional case}, considering the coordinate-based representation, yields:

\begin{equation}
    \tilde{x}_{2D} = \left ( (x_i-x_{mid})_{1\leq i \leq n},(y_i-y_{mid})_{1\leq i \leq n} \right ) \revisetwo{\in \mathbb{R}^{2n}},
\end{equation}
where $x_{mid}=\frac{1}{n}\sum_{i=1}^n x_i$ and $y_{mid}=\frac{1}{n}\sum_{i=1}^n y_i$.

The distance-based features yields the following feature vector:

\begin{equation}
    x_{2D} = \left( \Delta y,2 \Delta y,\Delta x,\sqrt{\Delta x ^2+\Delta y^2},\sqrt{\Delta x^2+(2 \Delta y)^2},2\Delta x,\sqrt{(2\Delta x)^2+\Delta y^2},2\sqrt{\Delta x^2+ \Delta y^2} \right).
\end{equation}
Using the same procedure as described for the 1D case, we consider the inverse of the distances, center and normalise all data points.

The structure of the neural network for $n=9$ and distance-based features in two dimensions has the following parameters: $p=4$, $d_0=8$ , $d_1=16$, $d_2=8$, $d_3=3$, $d_4=1$, with ReLU activation function for $p=1,2,3$ and linear activation for $p=4$. In 2D, we consider a few test cases and when the condition number of interpolation matrix is near $10^{12}$, we terminate training of the \reviseone{NN}.

\reviseone{\subsection{Integration of the data-driven method with interpolation and RBF-FD tasks}}

\reviseone{
In this section, we demonstrate how we use our data-driven model, that approximates \eqref{eq:approx function}, in practice. In particular, a sketch of the algorithm for the interpolation in an interval $[a,b]$ using $n$-points is detailed in algorithm \ref{alg:interpolation}. If we wish to interpolate on a set of points $N$ larger than $n$, we take subsets of $n$ points and evaluate the algorithm describe, with at least one point overlap to guarantee continuity of the interpolation.}

\begin{algorithm}
\reviseone{
\textbf{Input:} A target function $g$ to be interpolated and a set of ordered $n$ nodes $\vec{x} = (x_1,...,x_n)$ in some interval $[a,b]$ (or a set of data-points $\{x_i,g(x_i)\}_{i=1}^n$)\\
\textbf{Output:} An interpolation of $g$ in $n$ points in $[a,b]$\\
\SetAlgoLined
- Apply data transform as in \eqref{eq:transform}, $\vec{x}' = \mbox{transform}(\vec{x})$ \\
- Compute shape parameter $\mathcal{H}(\vec{x}') = \varepsilon$.\\
- Solve \eqref{eq:approximation_prob}:
\[
M(\vec{x},\varepsilon) \begin{bmatrix}
 \lambda\\
 \gamma
 \end{bmatrix}=
 \begin{bmatrix}
 \mathbf{f}\\
 0
 \end{bmatrix},
 \] \\
- Return approximation in interval $[a,b]$ using $n$-points: 
\[ g_{approx}(x) =  \sum_{i=1}^n \lambda_{i}\phi(||\mathbf{x}-\mathbf{x}_i||) +\sum_{k=1}^m \gamma_{k}\bar{p}_k(\mathbf{x}), \quad x\in [a,b]
\] 
\caption{RBF Interpolation using an adaptive shape parameter}
\label{alg:interpolation}}
\end{algorithm}

\reviseone{
In algorithm \ref{alg:fd}, we denote the integration of our shape parameter estimator in the context of a RBF-FD code. In particular, it is when assembling the PDE operator matrix \eqref{eq:diff_operator} that we require the shape parameters. }

\begin{algorithm}
\reviseone{
\textbf{Input:} Set of nodal points $\{x_i\}_{i=1}^N$ and set of function values $\{f(x_i)\}_{i=1}^N$.\\
\textbf{Output:} PDE operator matrix  \eqref{eq:diff_operator} \\
\SetAlgoLined
\caption{PDE operator with an adaptive shape parameter} 
\For{$x \in \{x_1,...,x_N\}$}{
- Find n-nearest neighbours of $x$, $N_x$ \\
- Compute shape parameter $\varepsilon = \mathcal{H}(N_x)$\\
- Generate interpolation matrix using $\varepsilon$
}
- Assemble differential operator $\mathcal{L}$
\label{alg:fd}}

\end{algorithm}

\section{Numerical results}
\label{se_results}

In this section, several problems are considered to demonstrate the accuracy of  using distance-based features in one and two  dimensions. We also present the results obtained with coordinate-based features in one dimension in Appendix \ref{appendix:coordinate_1D}. 

The MLP is trained with an initial learning rate of $\eta=10^{-4}$ and the regularization parameter is set to $\beta = 10^{-5}$. In the computations, the
oversampling parameter $r$ is set to 4 and \revisetwo{we add constant polynomials to the interpolation, i.e. $m=1$.}

The errors reported are using the $L_1$-error  between the
approximation and exact solutions:
\[ L_1\text{-error} = \frac{1}{N}\sum_{i=1}^N|u_e(\textbf{x}_i)-u_a(\textbf{x}_i) |\]
where $u_e$ and $u_a$ are the exact and approximate solutions, respectively.
\subsection{Numerical results for 1D test cases}
Following the method in section \ref{sec:model_data_driven}, we  train a neural network to predict the shape parameter $\varepsilon$. We consider two different kernels: the IMQRBF and the GARBF. We also vary the number of points considered to build the interpolation as $n = 5, 7, 10$.

\begin{figure}[H]
\begin{center}
\subfigure[Condition number]
{\includegraphics[width=0.45\textwidth]{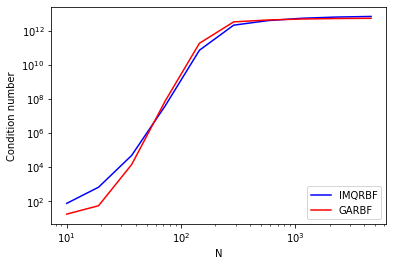}
}
\subfigure[Shape parameter]
{\includegraphics[width=0.45\textwidth]{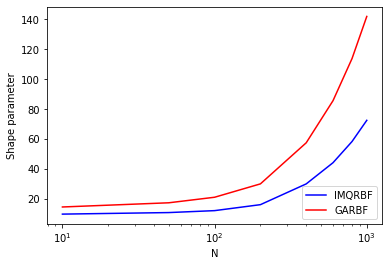}
}
\end{center}
	\caption{Condition number and shape parameter for different values of $N$ for IMQRBF and GARBF.}
	\label{fig:range_IMQRBF_GARBF}
\end{figure}

Fig. \ref{fig:range_IMQRBF_GARBF} demonstrates the performance of the trained \reviseone{NN}, using distance-based features, for $n=10$ stencil with IMQRBF and GARBF kernels. In the left panel, we show the average condition number of the interpolation matrix obtained as the distance between the nodal points diminishes (this occurs as we increase the number of points $N$ considered in the interval $[0,1]$). In the right panel, we show the predicted shape parameter again as the distance between nodal points decreases. In particular, we see that the predicted $\varepsilon$ increases as the distances decrease, as expected. 
Furthermore, we observe a similar behaviour when we consider $n=5$ and $7$.

In the next two sections, we will test this \reviseone{NN} on two specific tasks: \reviseone{function} interpolation (section \ref{sec:1d_interpolation}) and for solving a time-dependent problem (\ref{sec:1d_time_dep}).

\subsubsection{Interpolation}
\label{sec:1d_interpolation}
Two functions are used to examine the numerical accuracy of the present an \reviseone{adaptive} shape parameter
strategy in one-dimensional interpolation. We compare the proposed method with the constant shape parameter \reviseone{approach}, Hardy \cite{hardy}, Franke \cite{franke} and modified Franke strategies. For all the  interpolation tasks, the domain is discretized with the  equidistant and non-equidistant centers and the refined mesh is obtained by adding the midpoints in each step. We then consider sets of $n$ points to construct the interpolation, overlapping the boundary node between clusters.

The function $f_1$ is 
\[
f_1(x)=e^{\sin(\pi x)} \quad x  \in [0,1].    
\]

\begin{figure}[H]
	\centering
\includegraphics[width=0.45\textwidth]{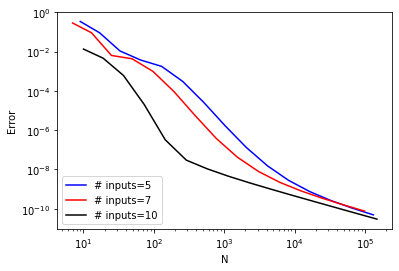}
	\caption{Convergence plot for the function $f_1$ using different input sizes for a neural network with IMQRBF.}
	\label{fig:f1_number_inputs}
\end{figure}

We first consider the IMQRBF kernel. Fig. \ref{fig:f1_number_inputs} shows the convergence plot for varying \reviseone{NN} input sizes for the function $f_1$ for equidistant points. Therefore, we can observe that an appropriate choice of  points is $n=10$.

\begin{figure}[H]
\begin{center}
\subfigure[equidistant]
{\includegraphics[width=0.45\textwidth]{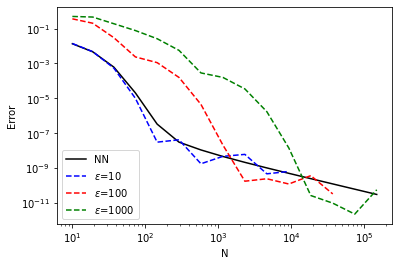}}
\subfigure[non-equidistant]
{\includegraphics[width=0.45\textwidth]{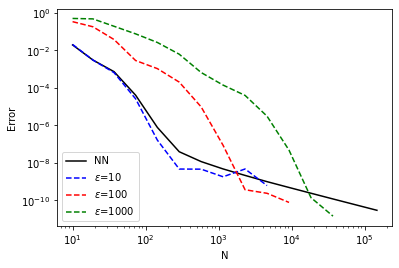}}
\end{center}
	\caption{Convergence plot for the function $f_1$ using a neural network strategy and constant shape parameters in 1D with IMQRBF.}
	\label{fig:f1}
\end{figure}

Fig. \ref{fig:f1} shows the  convergence plot for the function $f_1$ using shape parameters obtained by the \reviseone{NN} strategy and the constant shape parameters $\varepsilon = 10, 100, 1000$ as a function of the number of centers for equidistant and non-equidistant points. Note that as $N$ increases, when using \reviseone{a}  fixed $\varepsilon$, we run into the issue that the matrix becomes too ill-conditioned and we are unable to produce a suitable interpolation function, this is why, for example, for $\varepsilon=10$, we can only reach $N\approx 10^4$.

\begin{figure}[H]
\begin{center}
\subfigure[equidistant]
{\includegraphics[width=0.45\textwidth]{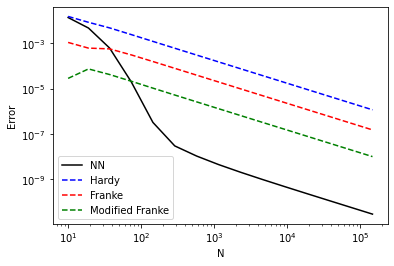}}
\subfigure[non-equidistant]
{\includegraphics[width=0.45\textwidth]{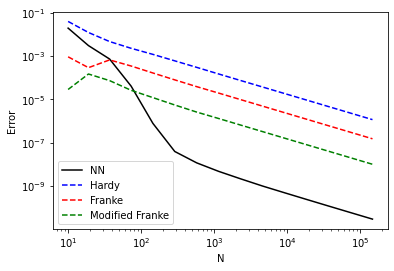}}
\end{center}
	\caption{Convergence plot for the function $f_1$ using a neural network, Hardy \cite{hardy}, Franke \cite{franke} and modified Franke strategies in 1D with IMQRBF.}
	\label{fig:f1_HF}
\end{figure}

Fig. \ref{fig:f1_HF}  compares the  convergence when using shape parameters obtained by the \reviseone{NN} and the adaptive methods Hardy, Franke and modified Franke strategies, as a function of the number of centers for equidistant and non-equidistant points. While these approaches always produce interpolation functions for this problem, the RBFs with
shape parameters obtained by \reviseone{NN} provide \reviseone{much better interpolation functions}  to interpolate 
function $f_1$.

\begin{figure}[H]
	\centering
\includegraphics[width=0.45\textwidth]{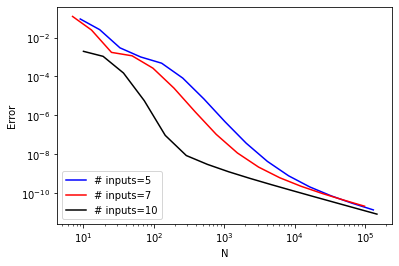}
	\caption{Convergence plot for the function $f_2$ using different input sizes for a neural network with IMQRBF.}
	\label{fig:f2_number_inputs}
\end{figure}

We now consider function $f_2$, the well known Runge function, 
\begin{equation*}
    f_2(x)=\frac{1}{1+16x^2},~x  \in [0,1].
\end{equation*}
Again, we focus first on the IMQRBF kernel. In Fig. \ref{fig:f2_number_inputs} we show the convergence plot for varying \reviseone{NN} input sizes for the function $f_2$ for equidistant points. The results look very similar to those shown in  \reviseone{Fig.} \ref{fig:f1_number_inputs}, and thus, we present the results using the $n=10$ \reviseone{for the NN}.

\begin{figure}[H]
\begin{center}
\subfigure[equidistant]
{\includegraphics[width=0.45\textwidth]{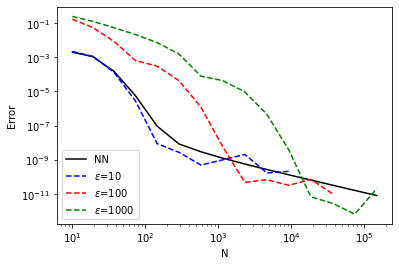}}
\subfigure[non-equidistant]
{\includegraphics[width=0.45\textwidth]{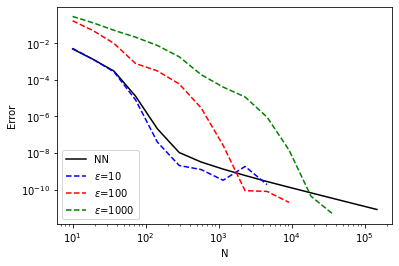}}
\end{center}
	\caption{Convergence plot for the function $f_2$ using a neural network strategy and constant shape parameters in 1D with IMQRBF.}
	\label{fig:f2}
\end{figure}

In Fig. \ref{fig:f2}, we show the  convergence plot using shape parameters obtained by the \reviseone{NN} strategy and constant shape parameters ($\varepsilon=10,100,1000$) as a function of the number of centers for equidistant and non-equidistant points. Again, we note that as $N$ increases, at some point the fixed $\varepsilon$ will fail as \reviseone{Fig.} \ref{fig:f1} produces interpolation matrices which are too ill-conditioned and we are unable to produce a suitable interpolation function.

\begin{figure}[H]
\begin{center}
\subfigure[equidistant]
{\includegraphics[width=0.45\textwidth]{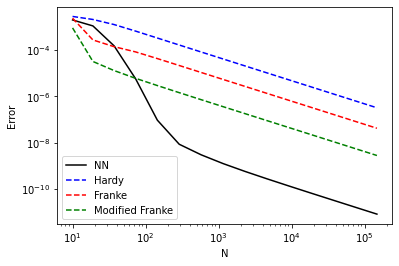}}
\subfigure[non-equidistant]
{\includegraphics[width=0.45\textwidth]{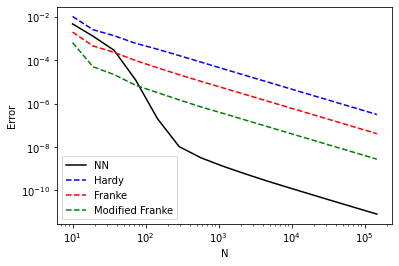}}
\end{center}
	\caption{Convergence plot for the function $f_2$ using a neural network, Hardy, Franke and modified Franke strategies in 1D with IMQRBF.}
	\label{fig:f2_HF}
\end{figure}

Fig. \ref{fig:f2_HF} shows the  convergence plot using shape parameters obtained by the \reviseone{NN}, Hardy, Franke and modified Franke strategies as a function of the number of centers for equidistant and non-equidistant points. Again, we observe much smaller error when using the \reviseone{NN}.

\begin{figure}[H]
\begin{center}
\subfigure[N=19]
{\includegraphics[width=0.45\textwidth]{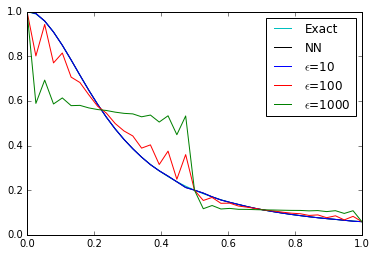}}
\subfigure[N=37]
{\includegraphics[width=0.45\textwidth]{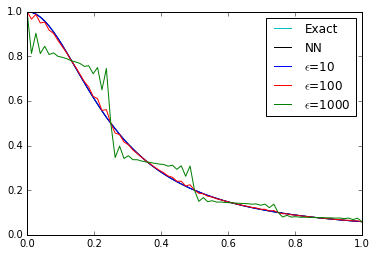}}
\end{center}
	\caption{Comparison between variable and constant shape parameter strategies for $f_2$ for $N=19$ and $N=37$ with IMQRBF.}
	\label{fig:sol_f2_uni_19_37}
\end{figure}

We also plot the resulting interpolation function for $N=19$ and $N=37$, comparing the variable and constant shape parameter strategies in Fig. \ref{fig:sol_f2_uni_19_37}.  We obtain similar results using the GARBF kernel for both interpolation tasks, which can be found in  Appendix \ref{1d_GARBF_interp}.


\begin{figure}[H]
\reviseone{
\begin{center}
\subfigure[equidistant]
{\includegraphics[width=0.45\textwidth]{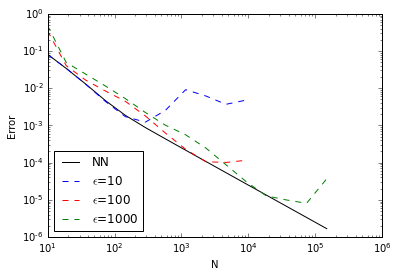}}
\subfigure[non-equidistant]
{\includegraphics[width=0.45\textwidth]{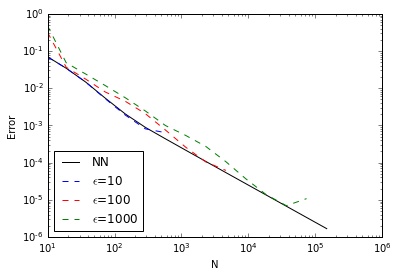}}
\end{center}
	\caption{\reviseone{Convergence plot for the function $f_3$ using a neural network strategy and constant shape parameters in 1D with IMQRBF.}}
	\label{fig:f5}
 }
\end{figure}

\begin{figure}[H]
\reviseone{
\begin{center}
\subfigure[equidistant]
{\includegraphics[width=0.45\textwidth]{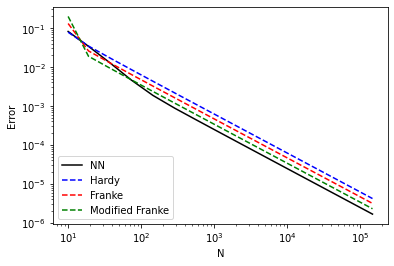}}
\subfigure[non-equidistant]
{\includegraphics[width=0.45\textwidth]{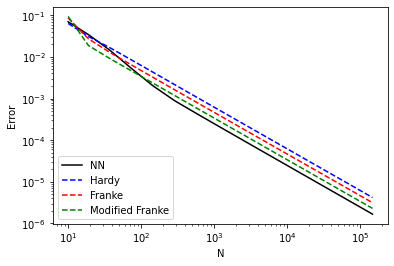}}
\end{center}
	\caption{\reviseone{Convergence plot for the function $f_3$ using a neural network, Hardy, Franke and modified Franke strategies in 1D with IMQRBF.}}
	\label{fig:f5_HF}
 }
\end{figure}

\begin{figure}[H]
\reviseone{
\begin{center}
\subfigure[N=19]
{\includegraphics[width=0.45\textwidth]{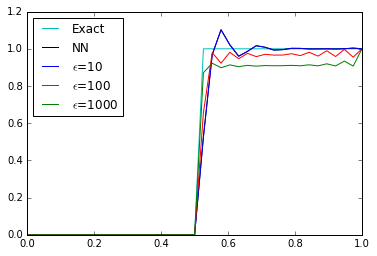}}
\subfigure[N=37]
{\includegraphics[width=0.45\textwidth]{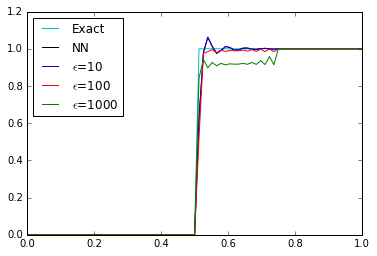}}
\end{center}
	\caption{\reviseone{Comparison between variable and constant shape parameter strategies for $f_3$ for $N=19$ and $N=37$ with IMQRBF.}}
	\label{fig:sol_f5_uni_19_37}
}
\end{figure}

\reviseone{The next considered test case in 1D for the interpolation is performed to assess the accuracy of our scheme on a piecewise constant problem as follows
\begin{equation*}
    f_3(x)=
\begin{cases}
1,&\text{if } x> 0.5, \\
0, &\text{if } x\leq 0.5,
\end{cases},~x  \in [0,1].
\end{equation*}
The convergence of the NN strategy and constant shape parameters $\varepsilon=10,100,1000$ is demonstrated in Fig \ref{fig:f5}. 
We observe
an overall good convergence  for our approach.   The results of the simulations comparing the NN,
Hardy, Franke and modified Franke schemes are illustrated in Fig. \ref{fig:f5_HF}. We
note that the error is smaller when using the NN, but the difference is smaller than the previous interpolation problems. The numerical solutions obtained with $N=19$ and $N=37$ are shown in Fig. \ref{fig:sol_f5_uni_19_37}.}

\subsubsection{Time-dependent problem}
\label{sec:1d_time_dep}
We have tested the NN strategy on a heat equation with two different initial conditions, and compare 
the RBF-FD method with the shape parameters obtained by the \reviseone{NN} strategy with finite difference method (FDM), finite element method (FEM) and RBF-FD method with constant shape parameters.  We  use BDF2 for time stepping in all cases. We use the same technique as in the previous section to obtain the refined mesh.
\\
The first  test case is the heat equation governed by 
\begin{equation}
\label{p1}
\frac{\partial u}{\partial t}=\frac{\partial ^2 u}{\partial x^2}, ~x  \in [0,1],  ~0 \leq t\leq 1,
\end{equation}
with the initial condition
\begin{equation}
\label{p2}
u(x,0)=-x^2+x,~x  \in (0,1) ,
\end{equation}
and the boundary conditions 
\begin{equation}
\label{p3}
u(0,t)=0,~ u(1,t)=0,~0\leq t\leq 1 .
\end{equation}
The exact solution of  \eqref{p1} with initial condition \eqref{p2} and boundary conditions \eqref{p3} is given by
$$u(x,t)=\sum_{n=1}^{\infty} \big(\frac{-4}{n^3 \pi^3}\big)\big((-1)^n-1\big)\sin\big(n\pi x\big)e^{-n^2 \pi ^2 t}.$$
We set $\Delta t=0.001$. The errors are shown in Table \ref{tab:heat_1d_1_uni} and Table \ref{tab:heat_1d_1_non_uni}  for equidistant and non-equidistant points, obtained with  \reviseone{the} IMQRBF kernel. Firstly, we show that \reviseone{in} our method, the error decreases consistently under mesh refinement for the equidistant points while it remains bounded for the non-equidistant points, attaining comparable performances to the FDM, FEM methods and RBF method with a fixed shape parameter $\varepsilon=10$. For the fixed shape parameter $\varepsilon=1$, we \reviseone{were} unable to attain a suitable numerical solution beyond $N=19$, and for $\varepsilon=100$, we observe a generally  approximation error. We obtain similar results  using a Gaussian kernel, shown in  Appendix \ref{1d_GARBF}. \revisetwo{The Table \ref{tab:time_heat_1d_1_uni} displays the execution time for equidistant points with the IMQRBF kernel. We see that with NN approach the execution time is of the same order of magnitude as when using a fixed shape parameter strategy. Furthermore, as the number of points $N$ increases, the cost of the NN method approaches the standard fixed shape parameter ones.}
\begin{table}[H]
  \centering
  \caption{The error using different schemes  for equidistant points with initial condition $\reviseone{u(x,0)=}-x^2+x$ with IMQRBF.} 
\begin{tabular}{|c|c c c c  c c |} 
\hline
N & FDM &FEM &  NN&  $\varepsilon=1$&$\varepsilon=10$ & $\varepsilon=100$\\
\hline
N=10& 1.5283e-04 & 1.5152e-04 &5.4703e-03 & 6.9435e-06& 5.9731e-03& 1.4708e-02 \\

N=19& 4.1141e-05 & 3.9723e-05 & 5.8909e-04& 8.6638e-07& 6.0268e-04&1.5342e-02  \\

N=37& 1.1152e-05 & 9.8825e-06 & 5.0922e-05& -&4.4987e-05 & 1.4919e-02 \\

N=73& 3.3997e-06 & 2.4375e-06 &4.4037e-06 & -& 2.9077e-06& 1.0560e-02\\

N=145& 1.4322e-06 &9.6827e-07  & 1.3377e-06& -& 9.0652e-07 & 1.1204e-03 \\
\hline
\end{tabular}
	\label{tab:heat_1d_1_uni}
\end{table}

\begin{table}[H]
  \centering
  \caption{The error using different schemes  for non-equidistant points with initial condition $\reviseone{u(x,0)=}-x^2+x$ with IMQRBF.} 
\begin{tabular}{|c|c c c c  c c |} 
\hline
N & FDM &FEM &  NN& $\varepsilon=1$& $\varepsilon=10$ & $\varepsilon=100$\\
\hline
N=10& 3.0018e-03 &3.1917e-03  &9.2538e-03 &2.9977e-03 &9.5229e-03 & 1.4744e-02 \\

N=19& 3.1257e-03 & 3.1769e-03 & 4.0141e-03&  3.1253e-03&4.0648e-03 &1.5521e-02 \\

N=37&3.2019e-03  & 3.2151e-03 &3.7890e-03 & -&3.7745e-03 & 1.5813e-02\\

N=73& 3.2435e-03 & 3.2469e-03 &3.3382e-03 & -& 3.2918e-03&1.4598e-02\\

N=145&3.2653e-03  & 3.2662e-03 & 3.2660e-03& - & 3.2652e-03&8.4020e-03\\
\hline
\end{tabular}
	\label{tab:heat_1d_1_non_uni}
\end{table}

\begin{table}[H]
  \centering\revisetwo{
  \caption{\revisetwo{The execution time in seconds using different schemes  for equidistant points with initial condition $u(x,0)=-x^2+x$ with IMQRBF.}}
\begin{tabular}{|c|c c  c c |} 
\hline
N &  NN&  $\varepsilon=1$&$\varepsilon=10$ & $\varepsilon=100$\\
\hline
N=10& 3.3922e+00 & 2.7765e+00&2.8063e+00 &  2.7961e+00 \\
\hline
N=19& 7.1963e+00 &5.9187e+00 &6.0301e+00 &  6.0007e+00 \\
\hline
N=37& 1.4507e+01 &- &1.2070e+01 & 1.2124e+01  \\
\hline
N=73& 3.0278e+01 & -& 3.2495e+01& 2.5168e+01  \\
\hline
N=145& 1.3790e+02 & -& 1.2648e+02& 1.2943e+02  \\
\hline
\end{tabular}
	\label{tab:time_heat_1d_1_uni}
}
\end{table}

As a second test, we consider again the  heat equation \eqref{p1} with the different initial condition
\begin{equation}
\label{p4}
u(x,0)=6\sin(\pi x),~x  \in (0,1)
\end{equation}
and the boundary conditions 
\begin{equation}
\label{p5}
u(0,t)=0,~ u(1,t)=0,~0\leq t\leq 1 .
\end{equation}
The exact solution is given by
$$u(x,t)=6\sin(\pi x)e^{-\pi^2 t}.$$
In Table \ref{tab:heat_1d_2_uni} and Table \ref{tab:heat_1d_2_non_uni} the  error\reviseone{s} are listed  for equidistant  and non-equidistant points  with IMQRBF, respectively. The differences between RBF-
FD with shape parameters obtained by neural network and RBF-FD with constant shape parameter are
not significant. \reviseone{The NN method} has similar performance  to $\varepsilon=10$ but a superior performance when compared   to $\varepsilon=1$ and $\varepsilon=100$. Also, the differences between RBF-
FD with shape parameters obtained by the \reviseone{NN}, FDM and FEM are not
significant for non-equidistant points. Moreover, we can see that the results of the new strategy, FDM and FEM are approximately the same for equidistant points as long as $N$ is not too low.

We obtain similar results when considering the GARBF kernel, which can be found in  Appendix \ref{1d_GARBF_timedependent}.
\begin{table}[H]
  \centering
  \caption{The error using different schemes  for equidistant points with initial condition $\reviseone{u(x,0)=}6\sin(\pi x)$ with IMQRBF.} 
\begin{tabular}{|c|c c c c  c c |} 
\hline
N & FDM &FEM &  NN&  $\varepsilon=1$& $\varepsilon=10$&  $\varepsilon=100$\\
\hline
N=10& 3.5303e-03 & 3.4600e-03 & 1.2696e-01 &  6.2511e-05& 1.3867e-01& 3.4160e-01 \\

N=19& 9.4294e-04 & 9.1108e-04 & 1.3518e-02 &  1.7241e-05&1.3833e-02 &3.5640e-01 \\

N=37& 2.5361e-04 & 2.2566e-04 & 1.1860e-03 &  -&1.0468e-03 &3.4678e-01\\

N=73& 7.5759e-05 & 5.4121e-05 & 1.0088e-04 & - & 6.5944e-05&2.4553e-01\\

N=145& 3.0626e-05 & 2.0276e-05 & 2.8705e-05 &  -&1.8624e-05 &2.6311e-02\\
\hline
\end{tabular}
\label{tab:heat_1d_2_uni}
\end{table}

\begin{table}[H]
  \centering
  \caption{The error using different schemes  for non-equidistant points with initial condition $\reviseone{u(x,0)=} 6\sin(\pi x)$ with IMQRBF.} 
\begin{tabular}{|c|c c c c  c c |} 
\hline
N & FDM &FEM &  NN&  $\varepsilon=1$& $\varepsilon=10$& $\varepsilon=100$\\
\hline
N=10& 6.9932e-02 & 7.4578e-02 & 2.1597e-01 & 6.9815e-02&2.2221e-01 & 3.4248e-01 \\

N=19& 7.2790e-02 & 7.4039e-02 & 9.3362e-02 & 7.2770e-02&9.4567e-02 & 3.6052e-01\\

N=37& 7.4552e-02 & 7.4875e-02 & 8.8268e-02 & -&8.7929e-02 &3.6743e-01\\

N=73& 7.5520e-02 & 7.5602e-02 & 7.7758e-02 &- & 7.6662e-02 &3.4017e-01\\

N=145& 7.6028e-02 & 7.6048e-02 & 7.6045e-02 & - & 7.6025e-02&1.9665e-01\\
\hline
\end{tabular}
\label{tab:heat_1d_2_non_uni}
\end{table}

\subsection{Numerical results for 2D test cases}

We generate the interpolation and evaluation point sets in two dimensions using regular Cartesian meshes. The results presented are using a modified version of a RBF-FD code \cite{igor-rbf-code, igor-rbf}.

\begin{figure}[H]
\begin{center}
\subfigure[Condition number]
{\includegraphics[width=0.45\textwidth]{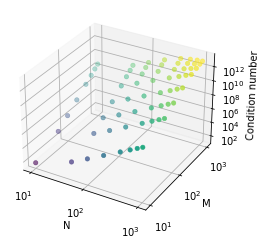}}
\subfigure[Shape parameter]
{\includegraphics[width=0.45\textwidth]{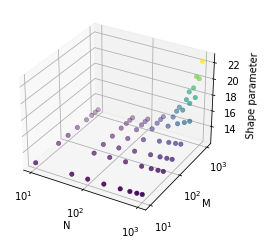}}
\end{center}
	\caption{Condition number and shape parameter for different values of $N$ and $M$ for IMQRBF.}
\label{fig:range_shape_2d}
\end{figure}

For the structured grid,  again, similarly to the results in one dimension, we note that the condition number remains well controlled for varying distances of points, and that the shape parameter increases as the distances between points decrease, as expected. We \reviseone{present} the condition number and shape parameter values for different values of $N$ and $M$  using the IMQRBF kernel in Fig.  \ref{fig:range_shape_2d}. We observe a similar behaviour when considering the GARBF kernel.

\subsubsection{Interpolation}
Our  numerical experiment of the present strategy, involves interpolating two functions in two dimensions. We approximate Franke’s function \cite{Fasshauer}
\begin{equation*}
\begin{split}
 f_4(x,y)&=\frac{3}{4} e^{-\big( \frac{(9x-2)^2+(9y-2)^2}{4}\big)}+\frac{3}{4} e^{-\big( \frac{(9x+1)^2}{49}+\frac{(9y+1)^2}{10}\big)}  
 \\&+\frac{1}{2} e^{-\big( \frac{(9x-7)^2+(9y-3)^2}{4}\big)}
 -\frac{1}{5} e^{-\big( (9x-4)^2+(9y-7)^2\big)}
 \end{split}
\end{equation*}
where $(x,y)\in [0,1]\times [0,1]$.  
The proposed NN method has been compared with a constant shape parameter strategy ($\varepsilon=1,5,10,100$).

We perform a convergence study for RBF-FD method using the GARBF kernel and we compare the results using the \reviseone{NN} and constant shape parameter strategies  ($\varepsilon=1,5,10,100$). The results are shown in Table \ref{tab:interpolation_2d_f3_structured_gaus}. We note that when using this kernel, there is a greater sensitivity to the choice of $\varepsilon$ and thus, we see a clear advantage to having an adaptive strategy \reviseone{for selecting} $\varepsilon$. 

\begin{table}[H]
  \centering
  \caption{The error using RBF-FD method using different strategies for interpolating $f_4$ for structured grid with GARBF.} 
\begin{tabular}{|c c |c c c  c c |} 
\hline
\# interpolation &  \#  evaluation&NN& $\varepsilon=1$&  $\varepsilon=5$& $\varepsilon=10$& $\varepsilon=100$\\
 points  &   points&& &  & & \\
\hline
10$\times$10& 20$\times$20 & 4.0537e-03 &4.0578e-03  &  3.9408e-03 & 7.4077e-03  &3.2520e-02 \\

20$\times$20& 40$\times$40 & 3.7430e-03 & 4.4900e-04 &  5.5438e-04 &2.3120e-03   & 1.2673e-02 \\

40$\times$40& 80$\times$80 & 5.4949e-04 & 6.8680e-05 &  6.7709e-05 & 3.0185e-04  &  5.1903e-03 \\

80$\times$80& 160$\times$160 & 7.0266e-05 & 7.9787e-03 & 8.2785e-06  & 3.7286e-05  &  9.3109e-04  \\

160$\times$160& 320$\times$320 & 9.2307e-06 & 1.6475e-01 & 4.8958e-06  & 4.6122e-06  &4.4077e-04
\\

320$\times$320& 640$\times$640 &1.3128e-06  & 1.5472e-01 & 1.2479e-03  &  4.8370e-06 &6.1821e-05 \\

\hline
\end{tabular}
	\label{tab:interpolation_2d_f3_structured_gaus}
\end{table}

\begin{figure}[H]
\begin{center}
\subfigure[NN]
{\includegraphics[width=0.45\textwidth]{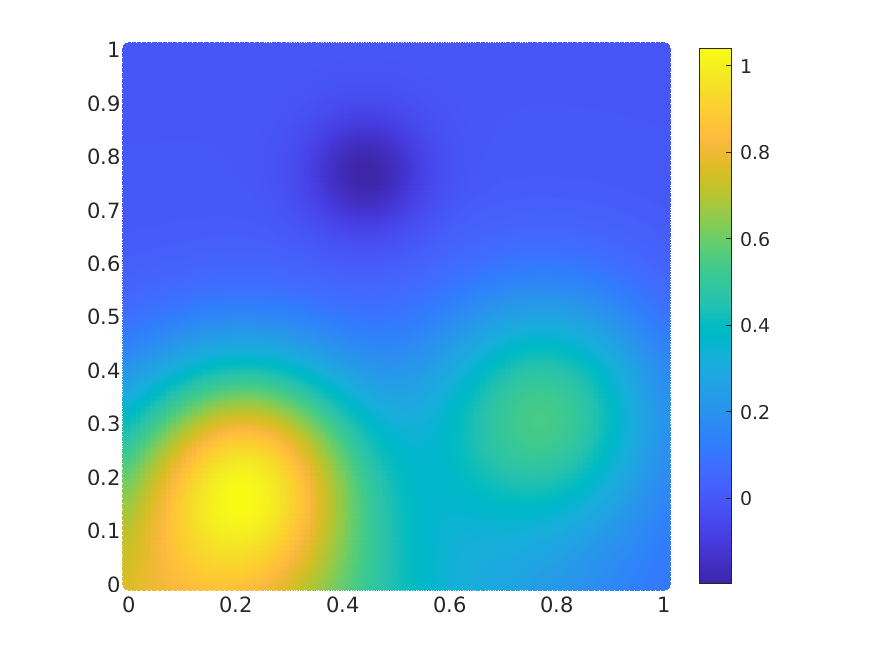}}
\subfigure[$\varepsilon=1$]
{\includegraphics[width=0.45\textwidth]{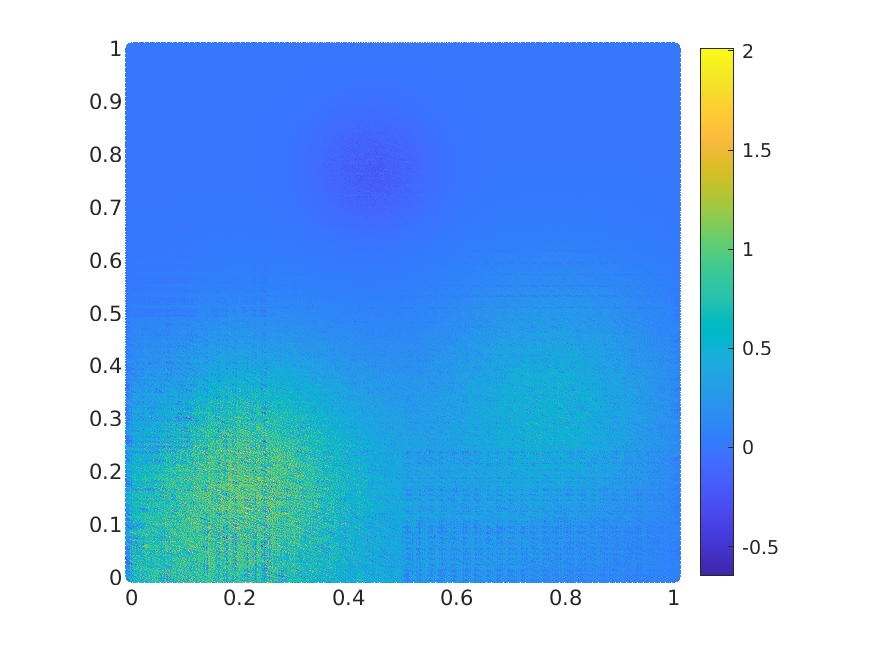}}
\subfigure[$\varepsilon=5$]
{\includegraphics[width=0.45\textwidth]{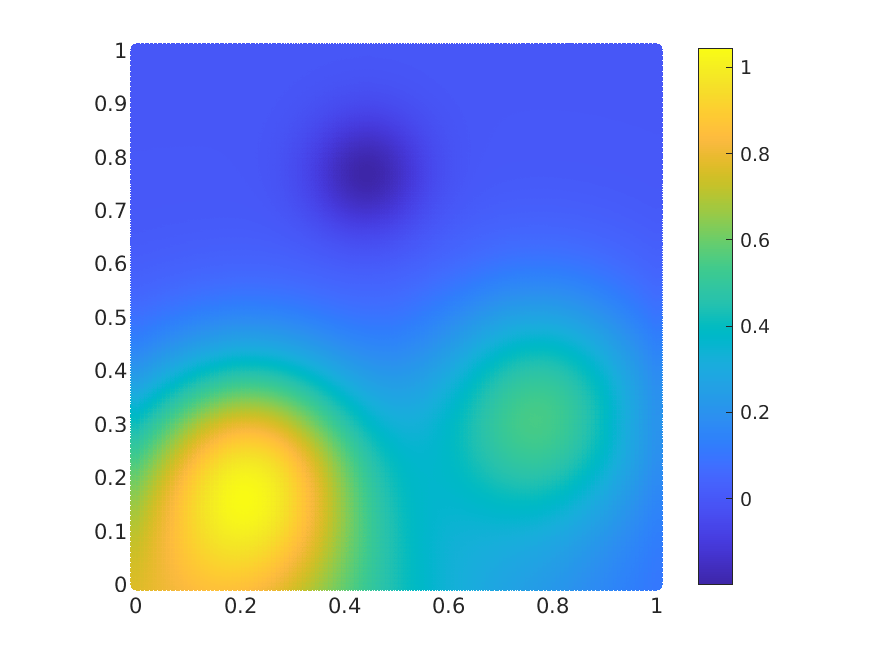}}
\subfigure[$\varepsilon=10$]
{\includegraphics[width=0.45\textwidth]{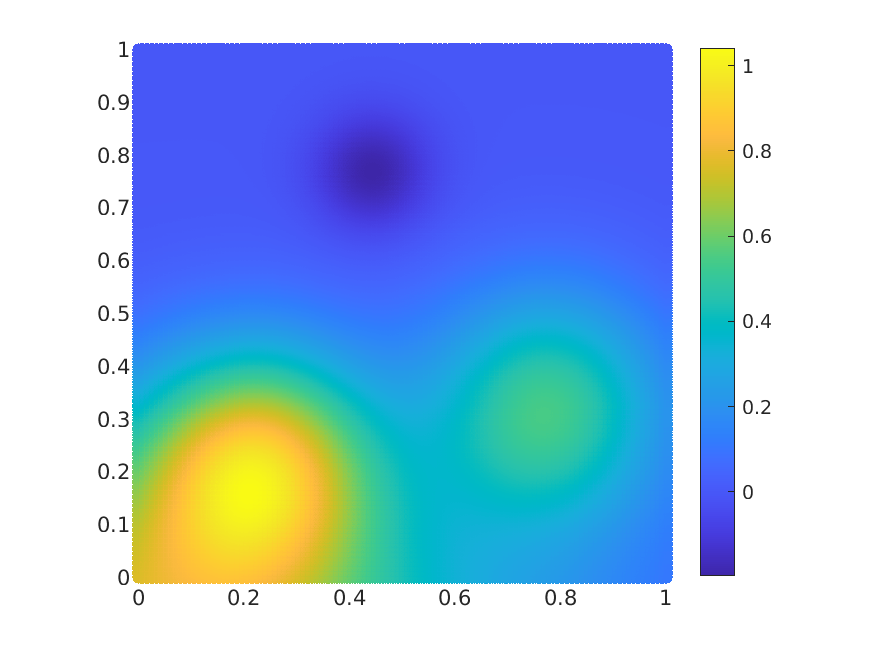}}
\end{center}
	\caption{Numerical solution for the function $f_4$ using $640\times 640$ evaluation points with IMQRBF.}
	\label{fig:f3}
\end{figure}

\begin{figure}[H]
\begin{center}
\subfigure[NN]
{\includegraphics[width=0.45\textwidth]{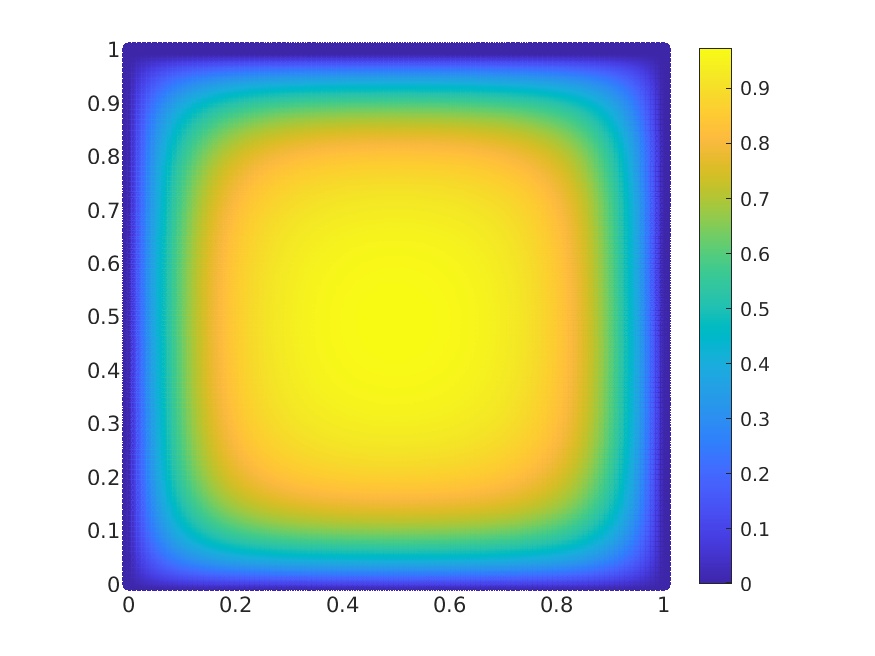}}
\subfigure[$\varepsilon=1$]
{\includegraphics[width=0.45\textwidth]{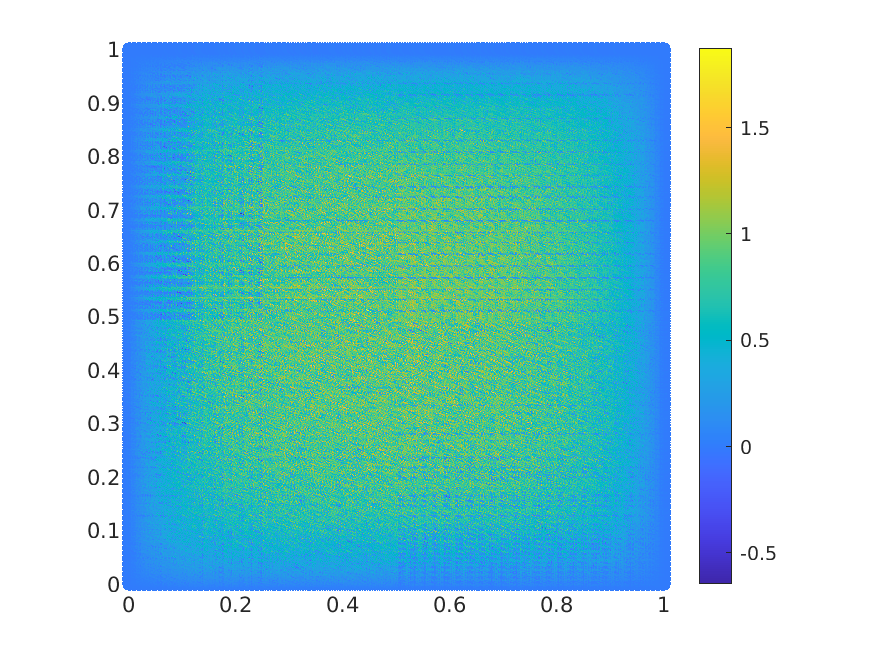}}
\subfigure[$\varepsilon=5$]
{\includegraphics[width=0.45\textwidth]{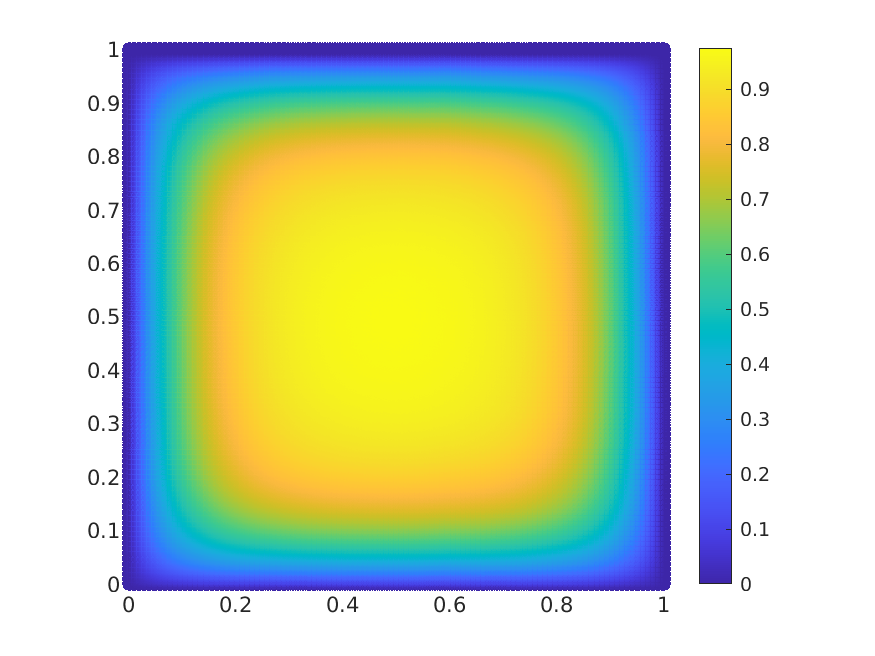}}
\subfigure[$\varepsilon=10$]
{\includegraphics[width=0.45\textwidth]{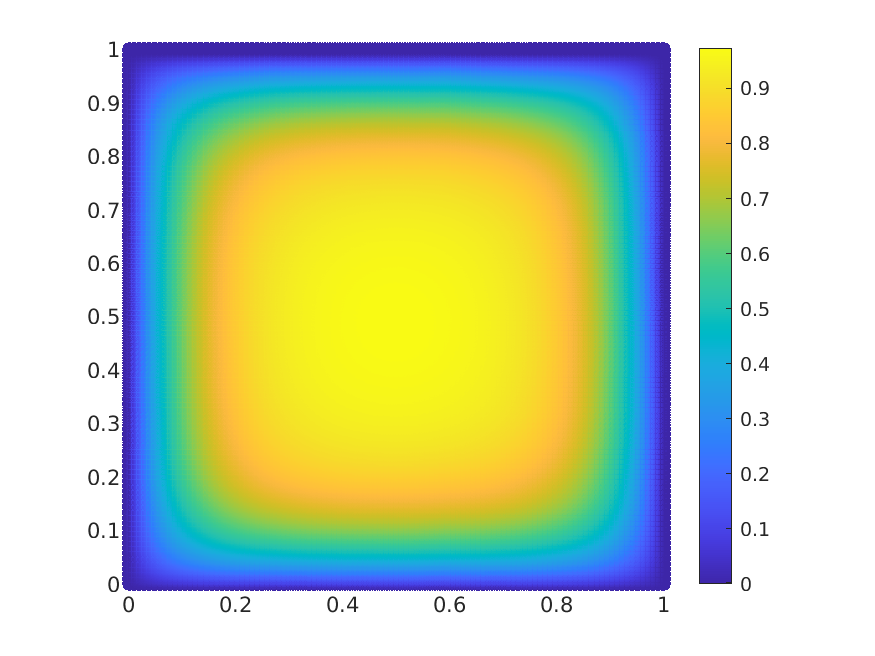}}
\end{center}
	\caption{Numerical solution for the function $f_5$ using $640\times 640$ evaluation points with IMQRBF.}
	\label{fig:f4}
\end{figure}

When considering the RBF-FD method using the IMQRBF kernel, we observe a similar error decay for NN and $\varepsilon=10$, as shown in Table \ref{tab:interpolation_2d_f3_structured}. Furthermore, we plot the numerical solution obtained by using NN and constant shape parameters with $640\times 640$ evaluation points in Fig. \ref{fig:f3}. We only report the results with $\varepsilon=1,5,10$, since the results for $\varepsilon=100$ \reviseone{are} of similar nature.
\begin{table}[ht]
  \centering
  \caption{The error using RBF-FD method using different strategies for interpolating $f_4$ for structured grid with IMQRBF.} 
\begin{tabular}{|c c |c c c  c c |} 
\hline
\# interpolation &  \#  evaluation&NN& $\varepsilon=1$&  $\varepsilon=5$& $\varepsilon=10$& $\varepsilon=100$\\
 points  &   points&& &  & & \\
\hline
10$\times$10& 20$\times$20 & 4.3784e-03 & 3.7428e-03 &3.3664e-03 &4.5244e-03 &2.3411e-02 \\

20$\times$20& 40$\times$40 & 2.3193e-03 & 4.2816e-04
 & 6.3583e-04&1.8377e-03 &5.5466e-03 \\

40$\times$40& 80$\times$80 & 5.1059e-04 & 5.0617e-05 &8.6087e-05 &3.3993e-04
 & 3.5434e-04 \\

80$\times$80& 160$\times$160 &  7.6687e-05& 1.0758e-03 &1.0795e-05 &4.6778e-05 & 7.2623e-04  \\

160$\times$160& 320$\times$320 & 1.0667e-05 & 3.6852e-02 &1.4960e-06 &5.9460e-06 & 3.0940e-04
\\

320$\times$320& 640$\times$640 & 1.5576e-06 & 7.3941e-02 &1.7185e-04
 &9.8128e-07 &6.3065e-05  \\

\hline
\end{tabular}
	\label{tab:interpolation_2d_f3_structured}
\end{table}
The next  test function is taken from  \cite{Fasshauer1,Gartland,Jankowska}. We consider the function 
\begin{equation*}
  f_5(x,y)= (1+ e^{-\frac{1}{\alpha}}-e^{-\frac{x}{\alpha}}-e^{\frac{(x-1)}{\alpha}})(1+ e^{-\frac{1}{\alpha}}-e^{-\frac{y}{\alpha}}-e^{\frac{(y-1)}{\alpha}}),~ (x,y) \in [0,1]\times [0,1].
\end{equation*}
As $\alpha>0$ decreases, the problem becomes more difficult.  
In Table \ref{tab:interpolation_2d_f4_structured_gaus} we present the
corresponding results using GARBF kernel for the same numbers of interpolation and evaluation points, whereas in Table \ref{tab:interpolation_2d_f4_structured} we present the
corresponding results using IMQRBFs. In terms of errors, they are very similar to the results \reviseone{of} the previous test case. The numerical solution obtained by using the \reviseone{NN method} and constant shape parameters with $640\times 640$ evaluation points for $\alpha=0.1$ is shown in Fig. \ref{fig:f4}.

\begin{table}[H]
  \centering
  \caption{The error using RBF-FD method using different strategies for interpolating $f_5$  for structured grid with GARBF.} 
\begin{tabular}{|c|c c |c c c  c c |} 
\hline
&\# interpolation  &  \#  evaluation&NN& $\varepsilon=1$&  $\varepsilon=5$& $\varepsilon=10$& $\varepsilon=100$\\
&points  &  points&& &  & & \\
\hline
$\alpha=1$& & & & & & &\\
&10$\times$10& 20$\times$20 & 1.9268e-04 & 1.5882e-05 &2.8331e-04 &4.6367e-04 & 1.4925e-03 \\

&20$\times$20& 40$\times$40 & 1.8951e-04 & 1.7094e-06 &3.4351e-05 &1.2386e-04 &  6.4021e-04 \\

&40$\times$40& 80$\times$80 & 2.7365e-05 & 1.5040e-06 & 4.0265e-06& 1.5767e-05&2.4397e-04  \\

&80$\times$80 & 160$\times$160  & 3.4516e-06 & 3.3276e-04 &4.8374e-07 &1.9183e-06 & 4.4864e-05 \\

&160$\times$160& 320$\times$320 & 4.4844e-07 & 7.0042e-03 & 1.9386e-07& 2.3545e-07&  2.0077e-05
\\

&320$\times$320& 640$\times$640 & 6.4760e-08 & 6.5279e-03 &5.1445e-05 &2.0242e-07 & 2.8334e-06 \\

\hline

$\alpha=0.1$& & & & & & &\\
&10$\times$10& 20$\times$20 & 2.6965e-02 & 1.1568e-02 &2.6705e-02 &3.7388e-02 & 9.7074e-02 \\

&20$\times$20& 40$\times$40 & 1.2069e-02 & 1.4181e-03 &3.3277e-03 &8.4496e-03 & 3.8001e-02  \\

&40$\times$40& 80$\times$80 &1.7349e-03  & 1.9894e-04 &3.8974e-04 &1.0688e-03 &  1.3912e-02 \\

&80$\times$80& 160$\times$160 &  2.1546e-04& 1.9756e-02 &4.6063e-05 &1.2804e-04 & 2.5090e-03 \\

&160$\times$160& 320$\times$320 & 2.7536e-05 &4.1583e-01  &1.3108e-05 & 1.5523e-05& 1.1248e-03
\\

&320$\times$320& 640$\times$640 &3.9773e-06  & 3.8633e-01 &3.0628e-03 &1.2252e-05 &  1.5898e-04
 \\

\hline
\end{tabular}
	\label{tab:interpolation_2d_f4_structured_gaus}
\end{table}

\begin{table}[H]
  \centering
  \caption{The error using RBF-FD method using different strategies for interpolating $f_5$  for structured grid with IMQRBF.} 
\begin{tabular}{|c|c c |c c c  c c |} 
\hline
&\# interpolation  &  \#  evaluation&NN& $\varepsilon=1$&  $\varepsilon=5$& $\varepsilon=10$& $\varepsilon=100$\\
&points  &  points&& &  & & \\
\hline
$\alpha=1$& & & & & & &\\
&10$\times$10& 20$\times$20 &  2.5032e-04& 3.1038e-05 &2.2328e-04 &2.8326e-04 & 1.0485e-03 \\
&20$\times$20& 40$\times$40 &  1.1963e-04& 3.3140e-06 &3.8237e-05 & 9.7860e-05&2.7322e-04 \\
&40$\times$40& 80$\times$80 &  2.5064e-05& 4.4151e-07 &4.9931e-06 &1.7203e-05 & 1.7237e-05\\
&80$\times$80 & 160$\times$160  & 3.7100e-06 & 4.0940e-05 &6.1578e-07 &2.3282e-06  &3.3546e-05\\

&160$\times$160& 320$\times$320 &  5.1081e-07& 1.5014e-03 &7.8305e-08 &2.9357e-07 &1.4084e-05
\\
&320$\times$320& 640$\times$640 &7.3937e-08
  & 2.7727e-03 &7.2542e-06 & 4.7251e-08& 2.8793e-06\\
\hline
$\alpha=0.1$& & & & & & &\\
&10$\times$10& 20$\times$20 &  2.7700e-02& 1.2855e-02 & 2.3944e-02& 2.8308e-02& 7.3125e-02\\

&20$\times$20& 40$\times$40 &  8.1087e-03& 1.5699e-03 &3.5878e-03 &6.9354e-03 &1.8655e-02 \\

&40$\times$40& 80$\times$80 &  1.5816e-03& 1.8545e-04 &4.4749e-04 & 1.1406e-03 &1.6996e-03 \\

&80$\times$80& 160$\times$160 &  2.2833e-04& 2.4572e-03 &5.3792e-05 &1.5055e-04 & 1.8506e-03\\
&160$\times$160& 320$\times$320 & 3.0922e-05 &8.9455e-02  &6.7846e-06 &1.8757e-05 &7.8172e-04
\\
&320$\times$320& 640$\times$640 & 4.4105e-06 &1.6994e-01  &4.3497e-04 &3.0998e-06 &1.6067e-04
\\
\hline
\end{tabular}
	\label{tab:interpolation_2d_f4_structured}
\end{table}
\subsubsection{Boundary value problem}
We use the two dimensional
linear elliptic boundary value problem
\begin{equation}
\label{poisson_eq_2d}
\begin{split}
u_{xx}+u_{yy}= {}&-4\pi^2\sin(2\pi x y)(x^2+y^2),~ (x,y) \in [0,1]\times [0,1],
\\
u(x,y)={}& \sin(2\pi x y),~(x,y)  \in \partial([0,1]\times [0,1])
\end{split}
\end{equation}
as a steady PDE test problem for the shape parameter strategies. The exact solution
of the boundary value problem is $\reviseone{u(x,y)=}\sin(2\pi x y)$. We use collocated RBF-FD for this problem.

When considering the GARBF kernel, we again observe a greater sensitivity of the method with respect to the chosen $\varepsilon$, namely, that not only there is no fixed  $\varepsilon$ \reviseone{value} that performs better than the adaptive $\varepsilon$ provided by the \reviseone{NN}, but also for poorly chosen $\varepsilon$ ($\varepsilon=1 , 5, 10$), the approximation eventually blows-up, as shown in Table \ref{tab:poisson_2d_1_structured_gaus}. \revisetwo{In Table \ref{tab:time_poisson_2d_1_structured_gaus} we have reported the execution time using different variable and constant shape parameter strategies and note that they are of the same order of magnitude, in particular as the number of interpolation and evaluation points increase.}

\begin{table}[H]
  \centering
  \caption{The error using RBF-FD method using different strategies for Poisson equation with exact solution $u(x,y)=\sin(2\pi xy)$ for structured grid with GARBF.} 
\begin{tabular}{|c| c c c  c c |} 
\hline
\# interpolation/ evaluation points  &  NN& $\varepsilon=1$&  $\varepsilon=5$& $\varepsilon=10$& $\varepsilon=100$\\
\hline
10$\times$10&1.9403e-01 & 2.5470e-03 & 4.9695e-02 &1.3311e-01 & 3.3161e-01  \\

20$\times$20& 8.4645e-02 &  6.8800e-04& 1.5314e-02 &5.5160e-02 &  3.6776e-01  \\

40$\times$40& 2.7915e-02 &  2.9304e-04& 4.0221e-03 & 1.6304e-02& 3.5368e-01  \\

80$\times$80& 7.6639e-03 &  1.7146e-02& 1.0174e-03 &4.2514e-03 & 2.0177e-01 \\

160$\times$160&  2.0625e-03& 4.8929e-01 & 2.7545e-04 &1.0735e-03 &  8.1534e-02 \\

320$\times$320& 5.9784e-04 & 5.4754e+03 & 9.1389e-03 &3.1306e-04 &  2.5303e-02\\

640$\times$640& 9.0457e-04 & 2.1826e+02 &  1.0231e+03
&2.8746e-02 &6.7143e-03
\\

\hline
\end{tabular}
	\label{tab:poisson_2d_1_structured_gaus}
\end{table}

\begin{table}[H]
  \centering\revisetwo{
  \caption{\revisetwo{The execution time in seconds using RBF-FD method using different strategies for Poisson equation with exact solution $u(x,y)=\sin(2\pi xy)$ for structured grid with GARBF.}} 
\begin{tabular}{|c| c c c  c c |} 
\hline
\# interpolation/ evaluation points  &  NN& $\varepsilon=1$&  $\varepsilon=5$& $\varepsilon=10$& $\varepsilon=100$\\
\hline
10$\times$10& 1.9884e-01 & 1.7182e-01 &  1.5800e-01
& 1.5345e-01& 1.5336e-01 \\
20$\times$20& 2.5946e-01 & 1.8130e-01 & 1.7281e-01 
&1.7415e-01 & 1.7230e-01 \\
40$\times$40& 4.6245e-01 & 2.6422e-01 &  2.5122e-01
& 2.6028e-01& 2.4827e-01 \\
80$\times$80& 1.3021e+00 & 1.0107e+00 &  5.6645e-01
& 5.7521e-01&5.8777e-01 \\
160$\times$160& 5.0198e+00 & 4.0641e+00 &  2.1892e+00
& 2.2481e+00& 2.2854e+00\\
320$\times$320& 2.7198e+01 & 2.3509e+01 &  1.6101e+01
&1.6482e+01 &1.6921e+01\\
640$\times$640& 3.4350e+02 & 3.2823e+02 &  3.2782e+02
& 3.0111e+02&3.0119e+02
\\
\hline
\end{tabular}\label{tab:time_poisson_2d_1_structured_gaus}
}
\end{table}

\begin{figure}[H]
\begin{center}
\subfigure[NN]
{\includegraphics[width=0.45\textwidth]{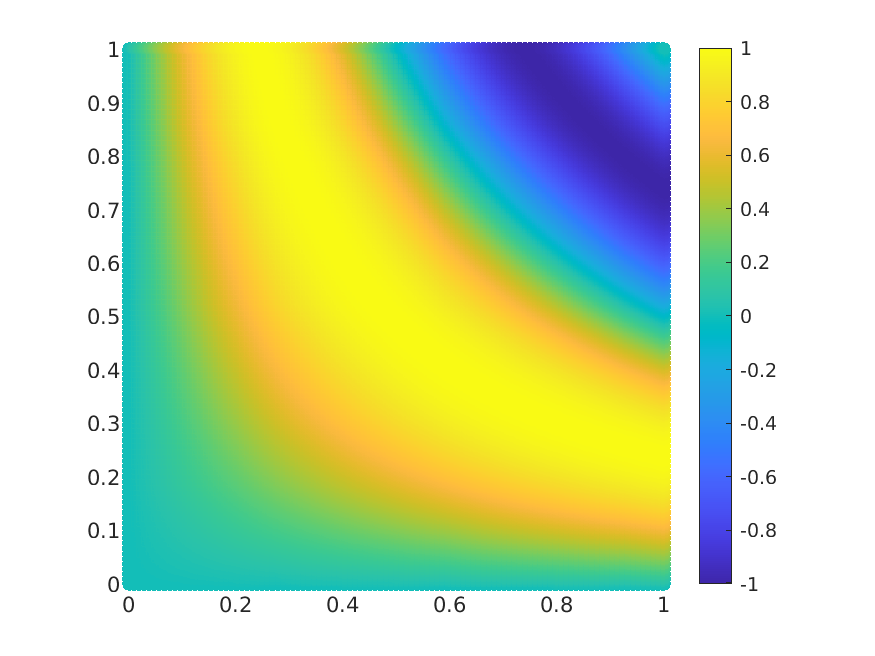}}
\subfigure[$\varepsilon=1$]
{\includegraphics[width=0.45\textwidth]{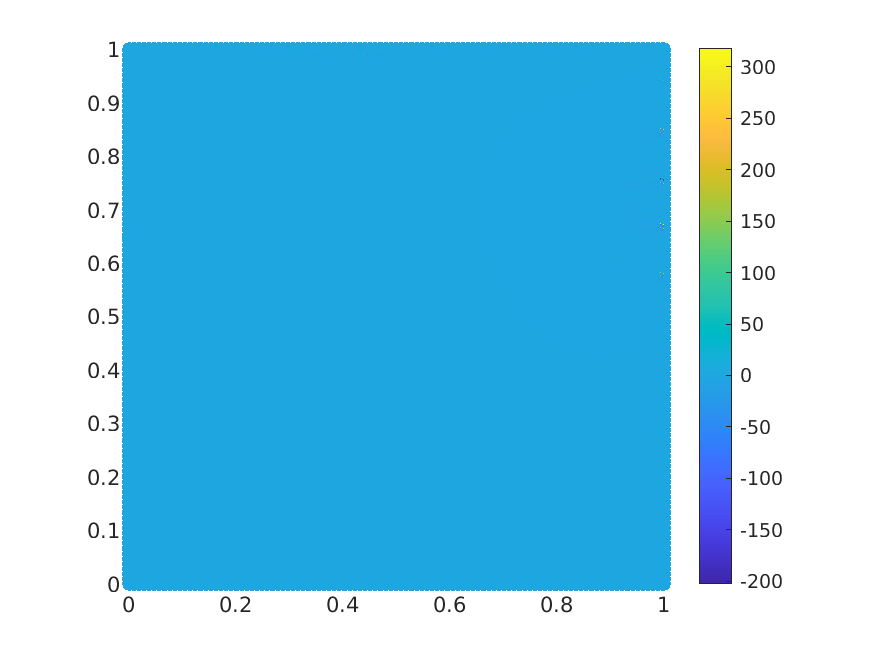}}
\subfigure[$\varepsilon=5$]
{\includegraphics[width=0.45\textwidth]{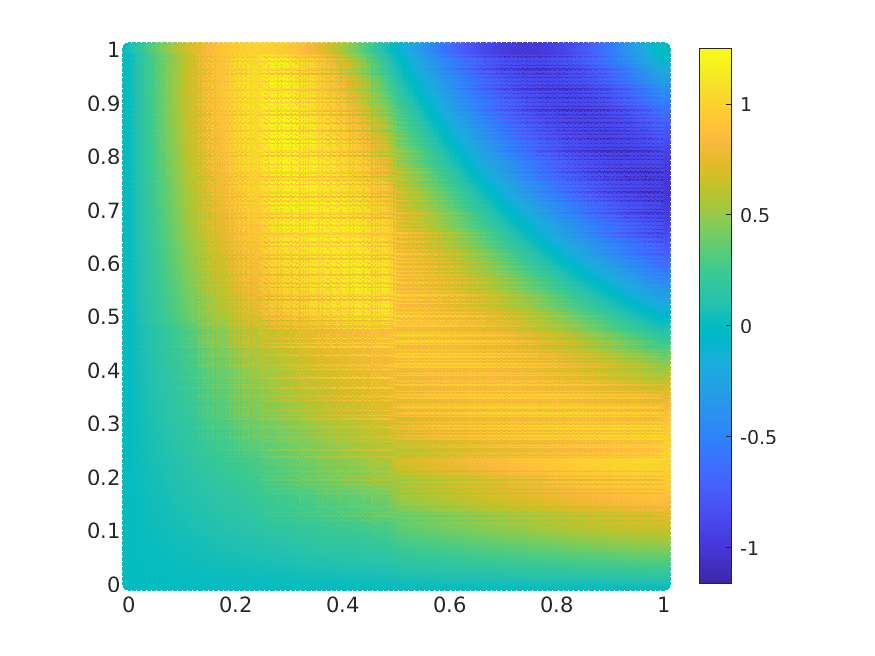}}
\subfigure[$\varepsilon=10$]
{\includegraphics[width=0.45\textwidth]{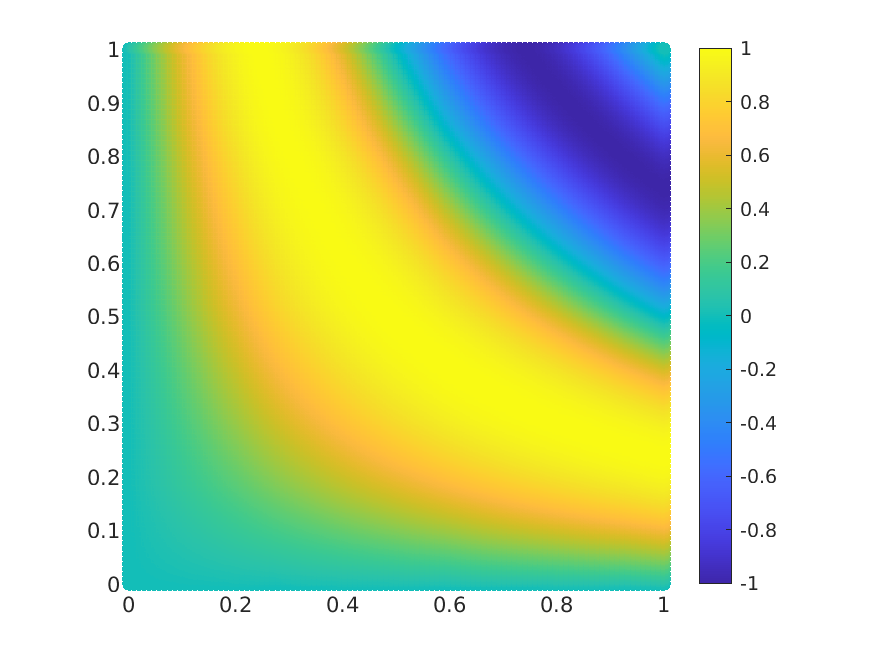}}
\end{center}
	\caption{Numerical solution for the Poisson equation using $640\times 640$ evaluation points with IMQRBF.}
	\label{fig:poisson}
\end{figure}

 The error when using the IMQRBF kernel with constant shape parameters ($\varepsilon=1, 5, 10, 100$) and variable shape parameter strategies are listed in Table \ref{tab:poisson_2d_1_structured}.  We see a clear advantage of the new method with respect
to the constant shape parameter strategy.  The errors blow-up in all the  cases except when using the NN. The numerical solution obtained by using NN and constant shape parameters with $640\times 640$ evaluation points is shown in the \ref{fig:poisson}. We  observe that the numerical solution is almost the same for NN and $\varepsilon=10$.
\begin{table}[H]
  \centering
  \caption{The error using RBF-FD method with different strategies for Poisson equation with exact solution $u(x,y)=\sin(2\pi xy)$ on a structured grid with IMQRBF.} 
\begin{tabular}{|c| c c c  c c |} 
\hline
\# interpolation/ evaluation points  &  NN& $\varepsilon=1$&  $\varepsilon=5$& $\varepsilon=10$& $\varepsilon=100$\\
\hline
10$\times$10& 1.8630e-01 &3.8201e-03  & 6.0522e-02 &1.4643e-01 & 3.2945e-01\\

20$\times$20& 8.7444e-02 &9.8806e-04  & 2.5272e-02 &6.5721e-01 &3.5997e-01 \\

40$\times$40&3.7178e-02  & 2.4926e-04 & 7.7340e-03 & 2.6150e-02&3.3186e-01\\

80$\times$80& 1.2454e-02 & 9.0683e-03 & 2.0447e-03 & 7.8927e-03&2.1501e-01\\

160$\times$160&3.6264e-03  &1.0339e-01  & 5.1784e-04 &2.0813e-03 &9.1097e-02\\

320$\times$320& 1.0744e-03 & 3.1287e+02 & 7.9589e-04 &5.2789e-04 &3.7253e-02\\

640$\times$640& 4.3698e-04 &7.0244e-01  & 6.4113e-02 & 8.9551e-04&1.1979e-02
\\
\hline
\end{tabular}
	\label{tab:poisson_2d_1_structured}
\end{table}
\section{Discussion and conclusion remarks}
\label{se_discussion}
We presented a methodology to train a \reviseone{NN} in an unsupervised fashion to predict a good shape parameter for given type of basis functions and size of stencil, taking as input location or relative distance of the interpolation nodes.

We studied different sizes of stencil $n=5,7,10$ and found that $n=10$ seems to lead to a better interpolation performance. This makes sense as we are using more interpolation points.

We also studied the different representation of the features. While coordinate-based features appeared more appealing and less costly for 2 and 3-dimension extensions, our numerical experiments demonstrate that using the relative distance between the interpolation points leads to a better performance (see Appendix \ref{appendix:coordinate_1D} for a coordinate-based features approach in 1D). Furthermore, the distance-based features are naturally shift invariant and can be made rotation invariant in higher dimensions.

We considered two types of radial basis kernels (IMQRBF and GARBF) plus a constant basis function. This methodology can be extended to different basis kernels without much effort, except retraining the network.

Still on the two considered radial basis kernels for constant shape parameter strategy, the performance of
GARBF is similar to \reviseone{that} of  IMQRBF, but in most  simulations the results  with \reviseone{the} IMQRBF are more robust  since
the blow-up happens later. In our experiments, we observe that the results for the \reviseone{NN} strategy with \reviseone{the} GARBF are much better than IMQRBF.

\revisetwo{In general, we have observed running times of the same magnitude when using our proposed method in comparison to a fixed shape parameter $\varepsilon$, with the gap becoming smaller as the number of points increases}. Furthermore, it is well known that RBF methods are, in general, very competitive when compared with traditional FDM/FEM \cite{zbMATH02009180,zbMATH06487347, zbMATH02082648}.

The presented methodology is very flexible in terms of \reviseone{the} size of the stencil, basis functions, desired matrix condition number range (in our work we set it to be between $10^{10}$ and $10^{12}$). The main challenges are the correct description of the input features and the dataset generation, for which we provide our methodology as well.

\subsection{Future work}
Overall, the different tasks used to test our \reviseone{NN} in 1D show that this methodology is robust a across large range of distributions in points, both when considering equidistant and non-equidistant points. In 2D, we presented a \reviseone{NN} that considers a simple grid-like distribution of points and estimates a suitable shape parameter $\varepsilon$.

Although we have shown that this proposed methodology is robust across different tasks, it requires a specific structure on the cloud points. In the future, we would like to consider generic clouds of points. This comes with a few difficulties, e.g.: if we consider the distance-based features, which seems like the more suitable choice due to easily guaranteeing shift and rotation invariance, it leads to a feature space of dimension $\binom{n}{2}$ for $n$ nodes, as opposed to a feature space of dimension $2n$ when considering coordinate-based features. One possibility would be to consider graph \reviseone{NNs} which are more suited for graph-like data, as we can interpret the cloud of points as graphs \cite{graphnns}.

Another issue is on how to best generate the training data. Ultimately, we are interested in integrating these methods with (local) RBF-FD type methods so there is often some structure in the stencil. In particular, it is important to have a good sampling of the interior of the domain, as well as \reviseone{on} the boundaries. This means that we have to identify strategies to generate a dataset that has a good covering of the input space. We can achieve this through a more detailed sampling strategy, balancing the number of points \reviseone{on} the boundary and in the interior.

\section*{Acknowledgments}
We take this opportunity to thank Fabian Mönkeberg for very helpful discussions.\\
\subsubsection*{Fundings}

FNM was supported by the Institute  of Mathematics at the University of Zürich.  MHV was supported by the Michigan Institute for Data Science research fellowship at the University of Michigan and the Department of Mathematics at the University of Michigan. PÖ was supported by the Gutenberg Research College, JGU Mainz.
\bibliographystyle{unsrt}
\bibliography{references}  

\appendix
\section{Results for coordinate-based features in one dimensional case}\label{appendix:coordinate_1D}
\subsection{One-dimensional interpolation}
\begin{figure}[H]
\begin{center}
\subfigure[equidistant]
{\includegraphics[width=0.45\textwidth]{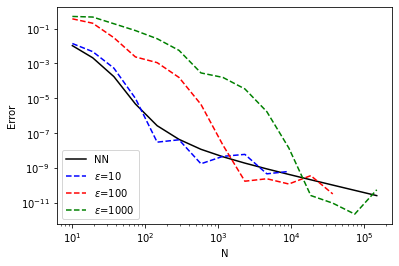}}
\subfigure[non-equidistant]
{\includegraphics[width=0.45\textwidth]{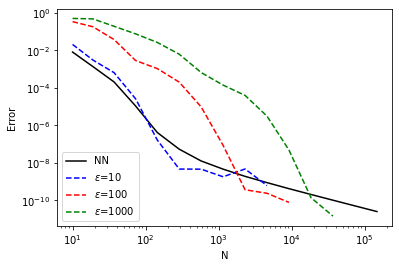}}
\end{center}
	\caption{Convergence plot for the function $f_1$ using a neural network strategy and constant shape parameters in 1D with IMQRBF.}
	\label{fig:f1_coord}
\end{figure}

\begin{figure}[H]
\begin{center}
\subfigure[equidistant]
{\includegraphics[width=0.45\textwidth]{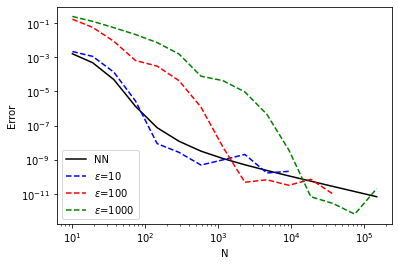}}
\subfigure[non-equidistant]
{\includegraphics[width=0.45\textwidth]{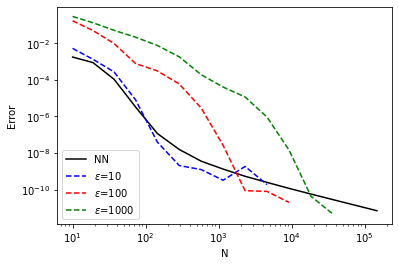}}
\end{center}
	\caption{Convergence plot for the function $f_2$ using a neural network strategy and constant shape parameters in 1D with IMQRBF.}
	\label{fig:f2_coord}
\end{figure}

Fig. \ref{fig:f1_coord} shows the  convergence plot for the function $f_1$ with IMQRBF using shape parameters obtained by the  \reviseone{NN} strategy based on coordinate-based features and the constant shape parameters as a function of the number of centers for equidistant and non-equidistant points.
Fig. \ref{fig:f2_coord} shows the  convergence plot for the function $f_2$ with IMQRBF using shape parameters obtained by the  \reviseone{NN} strategy based on coordinate-based features and the constant shape parameters as a function of the number of centers for equidistant and non-equidistant points. 
The results  for these two test cases are qualitatively similar to the ones obtained using distance-based features.

\subsection{One-dimensional time-dependent problem}
We compare the RBF-FD method with the shape parameters obtained by the \reviseone{NN} strategy based on coordinate-based features  with a FDM, a FEM and a RBF-FD method with constant shape parameters. In Table \ref{tab:heat_1d_1_uni_coord} and Table \ref{tab:heat_1d_1_non_uni_coord}, 
the  errors are listed  for  the heat equation \eqref{p1} with the  initial condition \eqref{p2}
for equidistant and non-equidistant points, respectively.
In Table \ref{tab:heat_1d_2_uni_coord} and Table \ref{tab:heat_1d_2_non_uni_coord} ,the  errors are listed  for  the heat equation \eqref{p1} with the  initial condition \eqref{p4} for equidistant and non-equidistant points, respectively. We  observe that we have the same results as for distance-based features for equidistant points, but the results for non-equidistant points \reviseone{are} better for distance-based features.
\begin{table}[H]
  \centering
  \caption{The error using different schemes  for equidistant points with initial condition $\reviseone{u(x,0)=}-x^2+x$ with IMQRBF.} 
\begin{tabular}{|c|c c c c  c c |} 
\hline
N & FDM &FEM &  NN&  $\varepsilon=1$&$\varepsilon=10$ & $\varepsilon=100$\\
\hline
N=10& 1.5283e-04 & 1.5152e-04 &1.9181e-03 & 6.9435e-06& 5.9731e-03& 1.4708e-02 \\

N=19& 4.1141e-05 & 3.9723e-05 &3.2476e-04 & 8.6638e-07& 6.0268e-04&1.5342e-02  \\

N=37& 1.1152e-05 & 9.8825e-06 &1.5181e-05 & -&4.4987e-05 & 1.4919e-02 \\

N=73& 3.3997e-06 & 2.4375e-06 &1.9688e-06 & -& 2.9077e-06& 1.0560e-02\\

N=145& 1.4322e-06 &9.6827e-07  &1.2384e-06 & -& 9.0652e-07 & 1.1204e-03 \\
\hline
\end{tabular}
	\label{tab:heat_1d_1_uni_coord}
\end{table}

\begin{table}[H]
  \centering
  \caption{The error using different schemes  for non-equidistant points with initial condition $\reviseone{u(x,0)=}-x^2+x$ with IMQRBF.} 
\begin{tabular}{|c|c c c c  c c |} 
\hline
N & FDM &FEM &  NN& $\varepsilon=1$& $\varepsilon=10$ & $\varepsilon=100$\\
\hline
N=10& 3.0018e-03 &3.1917e-03  &5.4642e-03 &2.9977e-03 &9.5229e-03 & 1.4744e-02 \\

N=19& 3.1257e-03 & 3.1769e-03 & 3.3695e-03&  3.1253e-03&4.0648e-03 &1.5521e-02 \\

N=37&3.2019e-03  & 3.2151e-03 &3.3639e-03 & -&3.7745e-03 & 1.5813e-02\\

N=73& 3.2435e-03 & 3.2469e-03 &3.2547e-03 &- & 3.2918e-03&1.4598e-02\\

N=145&3.2653e-03  & 3.2662e-03 & 3.2656e-03& - & 3.2652e-03&8.4020e-03\\
\hline
\end{tabular}
	\label{tab:heat_1d_1_non_uni_coord}
\end{table}

\begin{table}[H]
  \centering
  \caption{The error using different schemes  for equidistant points with initial condition $\reviseone{u(x,0)=}6\sin(\pi x)$ with IMQRBF.} 
\begin{tabular}{|c|c c c c  c c |} 
\hline
N & FDM &FEM &  NN&  $\varepsilon=1$& $\varepsilon=10$&  $\varepsilon=100$\\
\hline
N=10& 3.5303e-03 & 3.4600e-03 & 4.4323e-02 &  6.2511e-05& 1.3867e-01& 3.4160e-01 \\

N=19& 9.4294e-04 & 9.1108e-04 & 7.4256e-03 &  1.7241e-05&1.3833e-02 &3.5640e-01 \\

N=37& 2.5361e-04 & 2.2566e-04 & 3.5209e-04 &  -&1.0468e-03 &3.4678e-01\\

N=73& 7.5759e-05 & 5.4121e-05 & 4.3776e-05 &-  & 6.5944e-05&2.4553e-01\\

N=145& 3.0626e-05 & 2.0276e-05 & 2.6383e-05& - &1.8624e-05 &2.6311e-02\\
\hline
\end{tabular}
\label{tab:heat_1d_2_uni_coord}
\end{table}

\begin{table}[H]
  \centering
  \caption{The error using different schemes  for non-equidistant points with initial condition $\reviseone{u(x,0)=}6\sin(\pi x)$ with IMQRBF.} 
\begin{tabular}{|c|c c c c  c c |} 
\hline
N & FDM &FEM &  NN&  $\varepsilon=1$& $\varepsilon=10$& $\varepsilon=100$\\
\hline
N=10& 6.9932e-02 & 7.4578e-02 & 1.2766e-01 & 6.9815e-02&2.2221e-01 & 3.4248e-01 \\

N=19& 7.2790e-02 & 7.4039e-02 & 7.8271e-02 & 7.2770e-02&9.4567e-02 & 3.6052e-01\\

N=37& 7.4552e-02 & 7.4875e-02 &7.8321e-02 & -&8.7929e-02 &3.6743e-01\\

N=73& 7.5520e-02 & 7.5602e-02 & 7.5785e-02 &- & 7.6662e-02 &3.4017e-01\\

N=145& 7.6028e-02 & 7.6048e-02 & 7.6035e-02 & - & 7.6025e-02&1.9665e-01\\
\hline
\end{tabular}
\label{tab:heat_1d_2_non_uni_coord}
\end{table}

\section{Results for 1D GARBF}
\label{1d_GARBF}
\subsection{Interpolation}
\label{1d_GARBF_interp}
\begin{figure}[H]
\begin{center}
\subfigure[equidistant]
{\includegraphics[width=0.45\textwidth]{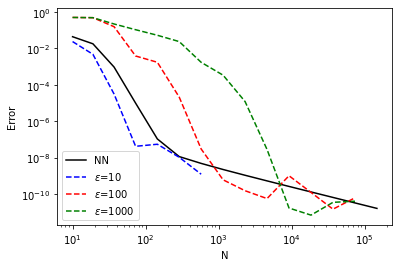}}
\subfigure[non-equidistant]
{\includegraphics[width=0.45\textwidth]{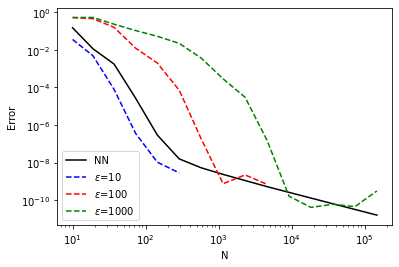}}
\end{center}
	\caption{Convergence plot for the function $f_1$ using a neural network strategy and constant shape parameters in 1D with GARBF.}
	\label{fig:f1_gaus}
\end{figure}

\begin{figure}[H]
\begin{center}
\subfigure[equidistant]
{\includegraphics[width=0.45\textwidth]{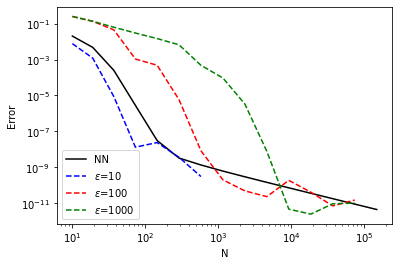}}
\subfigure[non-equidistant]
{\includegraphics[width=0.45\textwidth]{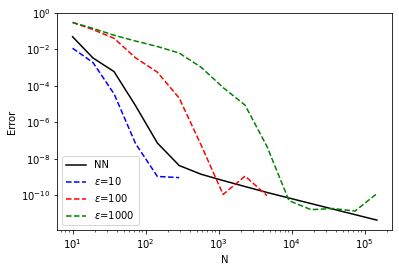}}
\end{center}
	\caption{Convergence plot for the function $f_2$ using a neural network strategy and constant shape parameters  in 1D with GARBF.}
	\label{fig:f2_gaus}
\end{figure}

 Fig. \ref{fig:f1_gaus}  shows the  convergence plot for the function $f_1$  with GARBF using shape parameters obtained by the \reviseone{NN} strategy and the constant shape parameters as a function of centers for equidistant and non-equidistant points. Fig. \ref{fig:f2_gaus}  shows the  convergence plot for the function $f_2$  with GARBF using shape parameters obtained by the \reviseone{NN} strategy and the constant shape parameters as a function of the number of centers for equidistant and non-equidistant points.  

\subsection{Time-dependent PDE}
\label{1d_GARBF_timedependent}
In this section we will focus on  \eqref{p1}, with  two different initial and boundary conditions using GARBFs. We consider the initial condition \eqref{p2} and boundary condition \eqref{p3}. \reviseone{The error using different schemes for equidistant and non-equidistant points are given in Table \ref{tab:heat_1d_1_uni_gaus} and Table \ref{tab:heat_1d_1_non_uni_gaus}, respectively.} And, the same problem with initial condition \eqref{p4} and boundary condition \eqref{p5}. \reviseone{The error using different schemes for equidistant and non-equidistant points listed in Table \ref{tab:heat_1d_2_uni_gaus} and Table \ref{tab:heat_1d_2_non_uni_gaus}, respectively.}
\begin{table}[H]
  \centering
  \caption{The error using different schemes  for equidistant points with initial condition $\reviseone{u(x,0)=}-x^2+x$ with GARBF.} 
\begin{tabular}{|c|c c c c  c c |} 
\hline
N & FDM &FEM &  NN& $\varepsilon=1$& $\varepsilon=10$& $\varepsilon=100$\\
\hline
N=10& 1.5283e-04 & 1.5152e-04 & 1.1940e-02& 1.2809e-06&3.7287e-03 & 1.4734e-02 \\

N=19& 4.1141e-05 & 3.9723e-05 & 1.7786e-03&  -& 7.5650e-04& 1.5486e-02 \\

N=37& 1.1152e-05 & 9.8825e-06 & 1.0877e-04& - &  5.8281e-06 &1.5440e-02 \\

N=73& 3.3997e-06 & 2.4375e-06 & 2.0636e-06&-  & 8.7901e-07  &1.0878e-02\\

N=145& 1.4322e-06 &9.6827e-07  & 9.5398e-07&- &- &5.2185e-03  \\
\hline
\end{tabular}
	\label{tab:heat_1d_1_uni_gaus}
\end{table}

\begin{table}[H]
  \centering
  \caption{The error using different schemes  for non-equidistant points with initial condition $\reviseone{u(x,0)=}-x^2+x$ with GARBF.} 
\begin{tabular}{|c|c c c c  c c |} 
\hline
N & FDM &FEM &  NN&  $\varepsilon=1$&$\varepsilon=10$&  $\varepsilon=100$\\
\hline
N=10& 3.0018e-03 &3.1917e-03  & 1.3410e-02& 2.9978e-03& 1.0159e-02& 1.4763e-02 \\

N=19& 3.1257e-03 & 3.1769e-03 & 6.5013e-03&  3.2367e-03& 4.0193e-03& 1.5594e-02 \\

N=37&3.2019e-03  & 3.2151e-03 & 4.4667e-03&  -& 3.2074e-03 &1.6058e-02 \\

N=73& 3.2435e-03 & 3.2469e-03 & 3.2711e-03&-  & 3.2434e-03 &1.5441e-02\\

N=145&3.2653e-03  & 3.2662e-03 & 3.2653e-03& - &  - & 1.1978e-02\\
\hline
\end{tabular}
\label{tab:heat_1d_1_non_uni_gaus}
\end{table}

\begin{table}[H]
  \centering
  \caption{The error using different schemes  for equidistant points with initial condition $\reviseone{u(x,0)=}6\sin(\pi x)$ with GARBF.} 
\begin{tabular}{|c|c c c c  c c |} 
\hline
N & FDM &FEM & NN&  $\varepsilon=1$& $\varepsilon=10$& $\varepsilon=100$\\
\hline
N=10& 3.5303e-03 & 3.4600e-03 &  2.7771e-01&  1.6327e-05&8.6403e-02 & 3.4220e-01 \\

N=19& 9.4294e-04 & 9.1108e-04 &  4.1225e-02&  -&1.7501e-02 &3.5974e-01 \\

N=37& 2.5361e-04 & 2.2566e-04 &  2.5212e-03&-  &1.3377e-04 &3.5886e-01\\

N=73& 7.5759e-05 & 5.4121e-05 &  4.6727e-05&- & 1.8014e-05&2.5313e-01\\

N=145& 3.0626e-05 & 2.0276e-05 &  1.9743e-05& - &- &1.2224e-01\\
\hline
\end{tabular}
\label{tab:heat_1d_2_uni_gaus}
\end{table}

\begin{table}[H]
  \centering
  \caption{The error using different schemes  for non-equidistant points with initial condition $\reviseone{u(x,0)=}6\sin(\pi x)$ with GARBF.} 
\begin{tabular}{|c|c c c c  c c |} 
\hline
N & FDM &FEM &  NN&  $\varepsilon=1$& $\varepsilon=10$& $\varepsilon=100$\\
\hline
N=10& 6.9932e-02 & 7.4578e-02 &  3.1237e-01& 6.9816e-02& 2.3717e-01&3.4291e-01  \\

N=19& 7.2790e-02 & 7.4039e-02 &  1.5131e-01& 7.2774e-02& 9.3198e-02& 3.6222e-01\\

N=37& 7.4552e-02 & 7.4875e-02 &  1.0438e-01& -&7.4683e-02 &3.7301e-01\\

N=73& 7.5520e-02 & 7.5602e-02 & 7.6176e-02 & - &7.5516e-02 &3.5991e-01\\

N=145& 7.6028e-02 & 7.6048e-02 &  7.6028e-02& - &- &2.8037e-01\\
\hline
\end{tabular}
\label{tab:heat_1d_2_non_uni_gaus}
\end{table}

\end{document}